\documentclass[11pt, oneside]{article} % use "amsart" for AMSLaTeX format

\usepackage[utf8]{inputenc}
\usepackage[T1]{fontenc}
\usepackage{amsmath}
\usepackage{amssymb} % note: it loads \usepackage{amsfonts}
\usepackage{mathrsfs} % for italic letters, use the command: \mathscr{}

\usepackage{amsthm} 
\theoremstyle{plain}
\newtheorem{theorem}{Theorem} % add [section] (or [chapter] in book class) to follow the section (chapter) number
\newtheorem{condition}[theorem]{Condition}
\newtheorem{example}[theorem]{Example}
\newtheorem{lemma}[theorem]{Lemma}
\newtheorem{proposition}[theorem]{Proposition}
\newtheorem{remark}[theorem]{Remark}

\usepackage{algorithm2e}
\SetKwInput{KwInputs}{Inputs}
\usepackage[noend]{algpseudocode}
\makeatletter
\def\BState{\State\hskip-\ALG@thistlm}
\makeatother

\usepackage{bm} % for isolated bold symbols within an equation
\usepackage{bbm} % for numbers in the \mathbb style
\RequirePackage[numbers]{natbib}
\RequirePackage[colorlinks,citecolor=blue,urlcolor=blue]{hyperref}

\newcommand{\N}{\mathbb N}

\newcommand{\R}{\mathbb R}

\usepackage{graphicx} % to insert images; use: \includegraphic{}
\usepackage{subcaption}

\usepackage{tikz} % to create plots with LaTeX
\usetikzlibrary{positioning ,shadows,matrix}
\usepackage[english]{babel} % to manage the language (you can specify more than one: e.g.[english,italian], and switch among them). It automatically handles the division in syllables

\usepackage{hyperref} % for hypertext references
\hypersetup{
    pdfpagemode={UseOutlines},
    bookmarksopen,
    pdfstartview={FitH},
    colorlinks,
    linkcolor={blue}, % color of the index
    citecolor={blue}, % color of bibliographic references
    urlcolor={blue} 
}

\usepackage{color} % to write in colors 

\usepackage[letterpaper,margin=1.25in]{geometry} % to set the margins
\def\iid{\overset{\textnormal{iid}}{\sim}} % i.i.d. symbol

\makeatletter  % cartesian product symbol
\@ifundefined{@currsizeindex} 
  {\RequirePackage{relsize}\let\dolarger\relsize} 
  {\def\dolarger#1{\larger[#1]}} 

\newcommand*\@@bigtimes[2]{\vphantom{\prod} 
  \vcenter{\hbox{\dolarger{4}$\m@th#1\mkern-2mu\times\mkern-2mu$}}} 
\newcommand*\bigtimes{\mathop{\mathpalette\@@bigtimes\relax}\displaylimits} 
\makeatother

\def\N{\mathbb{N}}\def\R{\mathbb{R}}\def\Z{\mathbb{Z}}\def\1{\mathbbm{1}}

\def\Acal{\mathcal{A}}\def\Bcal{\mathcal{B}}\def\Ccal{\mathcal{C}}\def\Dcal{\mathcal{D}}\def\Fcal{\mathcal{F}}\def\Hcal{\mathcal{H}}\def\Lcal{\mathcal{L}}\def\Ocal{\mathcal{O}}\def\Rcal{\mathcal{R}}\def\Tcal{\mathcal{T}}

\def\Gscr{\mathscr{G}}

\def\vol{\textnormal{vol}}

\title{Consistency of Bayesian inference with Gaussian process priors in an elliptic inverse problem}
\author{Matteo Giordano and Richard Nickl \\ \\ \textsc{University of Cambridge}}
\begin{document}

\maketitle

\abstract{
For $\Ocal$ a bounded domain in $\R^d$ and a given smooth function $g:\Ocal\to\R$, we consider the statistical nonlinear  inverse problem of recovering the conductivity $f>0$ in the divergence form equation
$$
	\nabla\cdot(f\nabla u)=g\ \textrm{on}\ \Ocal, \quad
	u=0\ \textrm{on}\ \partial\Ocal,
$$
from $N$ discrete noisy point evaluations of the solution $u=u_f$ on $\mathcal O$. We study the statistical performance of Bayesian nonparametric procedures based on a flexible class of Gaussian (or hierarchical Gaussian) process priors, whose implementation is feasible by MCMC methods. We show that, as the number $N$ of measurements increases, the resulting posterior distributions concentrate around the true parameter generating the data, and derive a convergence rate $N^{-\lambda}, \lambda>0,$ for the reconstruction error of the associated posterior means, in $L^2(\Ocal)$-distance.}

\bigskip


\tableofcontents

%------------------------------------------------------------------------------------------------------------------------
%------------------------------------------------------------------------------------------------------------------------
\section{Introduction}
\label{Sec:Intro}
%------------------------------------------------------------------------------------------------------------------------
%------------------------------------------------------------------------------------------------------------------------

%

Statistical inverse problems arise naturally in many applications in physics, imaging, tomography, and generally in engineering and throughout the sciences. A prototypical example involves a domain $\mathcal O \subset \mathbb R^d$, some function $f: \mathcal O \to \mathbb R$ of interest, and indirect measurements $G(f)$ of $f$, where $G$ is a given solution (or `forward') operator of some partial differential equation (PDE) governed by the unknown coefficient $f$. A natural statistical observational model postulates data
\begin{equation}\label{discrete}
Y_i = G(f)(X_i) + \sigma W_i,~~i=1, \dots, N,
\end{equation}
where the $X_i$'s are design points at which the PDE solution $G(f)$ is measured, and where the $W_i$'s are standard Gaussian noise variables scaled by a noise level $\sigma>0$. The aim is then to infer $f$ from the data $(Y_i, X_i)_{i=1}^N$. The study of problems of this type has a long history in applied mathematics, see the monographs \cite{EHN96, KNS08}, although explicit \textit{statistical} noise models have been considered only more recently \cite{KS04,BHM04,BHMR07, HP08}. Recent survey articles on the subject are \cite{BB18,AMOS19} where many more references can be found.

For many of the most natural PDEs -- such as the divergence form elliptic equation (\ref{Eq0}) considered below -- the resulting maps $G$ are \textit{non-linear} in $f$, and this poses various challenges: Among other things, the negative log-likelihood function associated to the model (\ref{discrete}), which equals the least squares criterion (see (\ref{Eq:JointLogLikelihood}) below for details), is then possibly \textit{non-convex}, and commonly used statistical algorithms (such as maximum likelihood estimators, Tikhonov regularisers or MAP estimates) defined as optimisers in $f$ of likelihood-based objective functions can not reliably be computed by standard convex optimisation techniques. While iterative optimisation methods (such as Landweber iteration) may overcome such challenges \cite{HNS95, Q00, KNS08, KSS09}, an attractive alternative methodology arises from the Bayesian approach to inverse problems advocated in an influential paper by Stuart \cite{S10}: One starts from a \textit{Gaussian process prior} $\Pi$ for the parameter $f$ or in fact, as is often necessary, for a suitable vector-space valued re-parameterisation $F$ of $f$. One then uses Bayes' theorem to infer the best posterior guess for $f$ given data $(Y_i, X_i)_{i=1}^N$. Posterior distributions and their expected values can be approximately computed via Markov Chain Monte Carlo (MCMC) methods (see, e.g., \cite{CRSW13, CMPS16, BGLFS17} and references therein) as soon as the forward map $G(\cdot)$ can be evaluated numerically, avoiding optimisation algorithms as well as the use of (potentially tedious, or non-existent) inversion formulas for $G^{-1}$; see Subsection \ref{Rem:Computation1} below for more discussion. The Bayesian approach has been particularly popular in application areas as it does not only deliver an estimator for the unknown parameter $f$ but simultaneously provides uncertainty quantification methodology for the recovery algorithm via the probability distribution of $f|(Y_i, X_i)_{i=1}^N$ (see, e.g., \cite{DS16}).  Conceptually related is the area of `probabilistic numerics' \cite{G19} in the noise-less case $\sigma=0$, with key ideas dating back to work by Diaconis \cite{D88}.

As successful as this approach may have proved to be in algorithmic practice, for the case when the forward map $G$ is non-linear we currently only have a limited understanding of the statistical validity of such Bayesian inversion methods. By validity we mean here \textit{statistical guarantees} for convergence of natural Bayesian estimators such as the posterior mean $\bar f = E^\Pi[f|(Y_i, X_i)_{i=1}^N]$ towards the ground truth $f_0$ generating the data. Without such guarantees, the interpretation of posterior based inferences remains vague: the randomness of the prior may have propagated into the posterior in a way that does not `wash out' even when very informative data is available (e.g., small noise variance and/or large sample size $N$), rendering  Bayesian methods potentially ambiguous for the purposes of valid statistical inference and uncertainty quantification.

	In the present article we attempt to advance our understanding of this problem area in the context of the following basic but representative example for a non-linear inverse problem: Let $g$ be a given smooth `source' function, and let $f: \mathcal O \to \mathbb R$ be a an unknown conductivity parameter determining solutions $u=u_f$ of the PDE
\begin{equation}
\label{Eq0}
	\begin{cases}
	\nabla\cdot(f\nabla u)=g & \textrm{on}\ \Ocal,\\
	u=0 & \textrm{on}\ \partial\Ocal,
	\end{cases}
\end{equation}
where we denote by $\nabla\cdot$ the divergence and by $\nabla$ the gradient operator, respectively. Under mild regularity conditions on $f$, and assuming that $f\ge K_{min}>0$ on $\Ocal$, standard elliptic theory implies that \eqref{Eq0} has a unique classical $C^2$-solution $G(f) \equiv u_f$. Identification of $f$ from an observed solution $u_f$ of this PDE has been considered in a large number of articles both in the applied mathematics and statistics communities -- we mention here \cite{R81, F83, HS85, KS85, A86, KL88, IK94, K01, S10, DS11, SS12, V13, DS16, BCDPW17, BGLFS17, NVW18, G19} and the many references therein.

	The main contributions of this article are as follows: We show that posterior means arising from a large class of Gaussian (or conditionally Gaussian) process priors for $f$ provide statistically consistent recovery (with explicit polynomial convergence rates as the number $N$ of measurements increases) of the unknown parameter $f$ in (\ref{Eq0}) from data in (\ref{discrete}).  While we employ the theory of posterior contraction from Bayesian non-parametric statistics \cite{vdVvZ08, vdVvZ09, GvdV17}, the non-linear nature of the problem at hand leads to substantial additional challenges arising from the fact that a) the Hellinger distance induced by the statistical experiment is not naturally compatible with relevant distances on the actual parameter $f$ and that b)  the `push-forward' prior induced on the information-theoretically relevant regression functions $G(f)$ is non-explicit (in particular, non-Gaussian) due to the non-linearity of the map $G$. Our proofs apply recent ideas from \cite{MNP19b} to the present elliptic situation. In the first step we show that the posterior distributions arising from the priors considered (optimally) solve the PDE-constrained regression problem of inferring $G(f)$ from data \eqref{discrete}. Such results can then be combined with a suitable ‘stability estimate' for the inverse map $G^{-1}$ to show that, for large sample size $N$, the posterior distributions concentrate around the true parameter generating the data at a convergence rate $N^{-\lambda}$ for some $\lambda>0$. We ultimately deduce the same rate of consistency for the posterior mean from quantitative uniform integrability arguments.

	The first results we obtain apply to a large class of `rescaled' Gaussian process priors similar to those considered in \cite{MNP19b}, addressing the need for additional a-priori regularisation of the posterior distribution in order to tame non-linear effects of the `forward map'. This rescaling of the Gaussian process depends on sample size $N$. From a non-asymptotic point of view this just reflects an adjustment of the covariance operator of the prior, but following \cite{D88} one may wonder whether a `fully Bayesian' solution of this non-linear inverse problem, based on a prior that does \textit{not} depend on $N$, is also possible. We show indeed that a hierarchical prior that randomises a finite truncation point in the Karhunen-Lo\'eve-type series expansion of the Gaussian base prior will also result in consistent recovery of the conductivity parameter $f$ in eq.~(\ref{Eq0}) from data (\ref{discrete}), at least if $f$ is smooth enough.

\smallskip

Let us finally discuss some related literature on statistical guarantees for Bayesian inversion: To the best of our knowledge, the only previous paper concerned with (frequentist) consistency of Bayesian inversion in the elliptic PDE (\ref{Eq0}) is by Vollmer \cite{V13}. The proofs in \cite{V13} share a similar general idea in that they rely on a preliminary treatment of the associated regression problem for $G(f)$, which is then combined with a suitable stability estimate for $G^{-1}$. However, the convergence rates obtained in \cite{V13} are only implicitly given and sub-optimal, also (unlike ours) for `prediction risk' in the PDE-constrained regression problem. Moreover, when specialised to the concrete non-linear elliptic problem (\ref{Eq0}) considered here, the results in Section 4 in \cite{V13} only hold for priors with bounded $C^\beta$-norms, such as `uniform wavelet type priors', similar to the ones used in \cite{NS17, N17, NS19} for different non-linear inverse problems. In contrast, our results hold for the more practical Gaussian process priors which are commonly used in applications, and which permit the use of tailor-made MCMC methodology -- such as the pCN algorithm discussed in Subsection \ref{Rem:Computation1} -- for computation.

The results obtained in \cite{NVW18} for the maximum a posteriori (MAP) estimates associated to the priors studied here are closely related to our findings in several ways. Ultimately the proof methods in \cite{NVW18} are, however, based on variational methods and hence entirely different from the Bayesian ideas underlying our results. Moreover, the MAP estimates in \cite{NVW18} are difficult to compute due to the lack of convexity of the forward map, whereas posterior means arising from Gaussian process priors admit explicit computational guarantees, see \cite{HSV14} and also Subsection \ref{Rem:Computation1} for more details.

It is  further of interest to compare our results to those recently obtained in \cite{AN19}, where the statistical version of the \textit{Cald\'eron problem} is studied. There the `Dirichlet-to-Neumann map' of solutions to the PDE (\ref{Eq0}) is observed, corrupted by appropriate Gaussian matrix noise. In this case, as only boundary measurements of $u_f$ at $\partial \Ocal$ are available, the statistical convergence rates are only of order $\log^{-\gamma} (N)$ for some $\gamma>0$ (as $N \to \infty$), whereas our results show that when interior measurements of $u_f$ are available throughout $\Ocal$, the recovery rates improve to $N^{-\lambda}$ for some $\lambda>0$.

There is of course a large literature on consistency of Bayesian \textit{linear} inverse problems with Gaussian priors, we only mention \cite{K11, R13, S13, KLS16, MNP17} and references therein. The non-linear case considered here is fundamentally more challenging  and cannot be treated by the techniques from these papers -- however, some of the general theory we develop in the appendix provides novel proof methods also for the linear setting. 
	
	\smallskip

This paper is structured as follows. Section \ref{main} contains all the main results for the inverse problem arising with the PDE model (\ref{Eq0}).  The proofs, which also include some theory for general non-linear inverse problems that is of independent interest, are given in Section \ref{Sec:Proofs} and Appendix \ref{App:GenInvProbl}. Finally, Appendix \ref{App:BoringThings} provides additional details on some facts used throughout the paper.

\section{Main results}\label{main}

\subsection{A statistical inverse problem with elliptic PDEs}
\label{Subsec:PrelimAndNotation}
%------------------------------------------------------------------------------------------------------------------------
%------------------------------------------------------------------------------------------------------------------------

\subsubsection{Main notation}
	Throughout the paper, $\Ocal\subset\R^d,\ d \in \mathbb N$, is a given nonempty open and bounded set with smooth boundary $\partial\Ocal$ and closure $\bar\Ocal$.

	The spaces of continuous functions defined on $\Ocal$ and $\bar\Ocal$ are respectively denoted $C(\Ocal)$ and $C(\bar\Ocal)$, and endowed with the supremum norm $\|\cdot\|_\infty$. For positive integers $\beta\in\N$, $C^\beta(\Ocal)$ is the space of $\beta$-times 
differentiable functions with uniformly continuous derivatives; for non-integer $\beta>0$, $C^\beta(\Ocal)$ is defined as
$$
	C^\beta(\Ocal)
	=
		\Bigg\{f\in C^{\lfloor\beta\rfloor}(\Ocal):\forall |i| = \lfloor \beta\rfloor,
		\sup_{x,y\in\Ocal,x\neq y}
		\frac{|D^i f(x)-D^i f(y)|}{|x-y|^{\beta-\lfloor \beta\rfloor}}<\infty\Bigg\},
$$
where $\lfloor \beta\rfloor$ denotes the largest integer less than or equal to $\beta$, and for any multi-index $i=(i_1,\dots,i_d)$, $D^i$ is the $i$-th partial differential operator. $C^\beta(\Ocal)$ is normed by
$$
	\|f\|_{C^\beta(\Ocal)}
	=
		\sum_{|i|\le \lfloor \beta\rfloor}\sup_{x\in\Ocal}|D^i f(x)|
		+\sum_{|i| = \lfloor 	\beta\rfloor} 
		\sup_{x,y\in\Ocal,\ x\neq y}\frac{|D^if(x)-D^if(y)|}{|x-y|^{\beta-\lfloor \beta\rfloor}},
$$
where the second summand is removed for integer $\beta$. We denote by $C^\infty(\Ocal)=\cap_{\beta}C^\beta(\Ocal)$ the set of smooth functions, and by $C^\infty_c(\Ocal)$ the subspace of elements in $C^\infty(\Ocal)$ with compact support contained in $\Ocal$.

Denote by $L^2(\mathcal O)$ the Hilbert space of square integrable functions on $\Ocal$, equipped with its usual inner product $\langle \cdot, \cdot \rangle_{L^2(\Ocal)}$. For  integer $\alpha\ge0$, the order-$\alpha$ Sobolev space on $\Ocal$ is the separable Hilbert space
$$
	H^\alpha(\Ocal)
	=
		\{f\in L^2(\Ocal) : \forall |i|\le\alpha, \ \exists\ D^if\in L^2(\Ocal)\},\ 
		\langle f,g\rangle_{H^\alpha(\Ocal)}
	=
		\sum_{|i|\le\alpha}\langle D^i f,D^i g\rangle_{L^2(\Ocal)}.
$$
For non-integer $\alpha\ge0$, $H^\alpha(\Ocal)$ can be defined by interpolation, see, e.g., \cite{Lions1972}. For any $\alpha\ge0$, $H^\alpha_c(\Ocal)$ will denote the completion of $C^\infty_c(\Ocal)$ with respect to the norm $\|\cdot\|_{H^\alpha(\Ocal)}$. Finally, if $K$ is a nonempty compact subset of $\Ocal$, we denote  by $ H^\alpha_K(\Ocal)$ the closed subspace of functions in $ H^\alpha(\Ocal)$ with support contained in $K$. Whenever there is no risk of confusion, we will omit the reference to the underlying domain $\Ocal$.

	Throughout, we use the symbols $\lesssim$ and $\gtrsim$ for inequalities holding up to a universal constant. Also, for two real sequences $(a_N)$ and $(b_N),$ we say that $a_N\simeq b_N$ if both $a_N\lesssim b_N$ and $b_N\lesssim a_N$ for all $N$ large enough. For a sequence of random variables $Z_N$ we write $Z_N = O_{\Pr}(a_N)$ if for all $\varepsilon>0$ there exists $M_\varepsilon<\infty$ such that for all $N$ large enough, $\Pr(|Z_N| \ge M_\varepsilon a_N) <\varepsilon$. Finally, we will denote by $\mathcal L(Z)$ the law of a random variable $Z$.

\subsubsection{Parameter spaces and link functions}
\label{Subsec:DivFormPDE}
%------------------------------------------------------------------------------------------------------------------------
%------------------------------------------------------------------------------------------------------------------------

	Let $g\in C^\infty(\Ocal)$ be an arbitrary source function, which will be regarded as fixed throughout. For $f\in C^\beta(\Ocal), \ \beta>1,$ consider the boundary value problem
\begin{equation}
\label{Eq:DivFormPDE}
	\begin{cases}
	\nabla\cdot(f\nabla u)=g & \textnormal{on}\ \Ocal,\\
	u=0 & \textnormal{on}\ \partial\Ocal.
	\end{cases}
\end{equation}
If we assume that $f\ge K_{min}>0$ on $\Ocal$, then standard elliptic theory (e.g., \cite{GT98}) implies that \eqref{Eq:DivFormPDE} has a classical solution $G(f)\equiv u_f \in C(\bar\Ocal)\cap C^{1+\beta}(\Ocal)$.

	We consider the following parameter space for $f$: For integer $\alpha>1+d/2, \ K_{min}\in (0,1)$, and denoting by $n=n(x)$ the outward pointing normal at $x\in\partial\Ocal$, let
\begin{equation}
\label{Eq:ParamSp}
	\Fcal_{\alpha,K_{min}}
	=
		\Bigg\{f\in H^\alpha(\Ocal): \textnormal{$\inf_{x \in \Ocal}f(x)>K_{min}$, 
		$f_{|\partial \Ocal}=1$, $\frac{\partial^j f}{\partial n^j}_{|\partial \Ocal}=0$
		for $1\le j\le\alpha-1$}\Bigg\}.
\end{equation}

	Our approach will be to place a prior probability measure on the unknown conductivity $f$ and base our inference on the posterior distribution of $f$ given noisy observations of $G(f)$, via Bayes' theorem. It is of interest to use \textit{Gaussian process priors}. Such probability measures are naturally supported in linear spaces (in our case $H^\alpha_c(\mathcal O)$) and we now introduce a bijective re-parametrisation so that the prior for $f$ is supported in the relevant parameter space $\mathcal F_{\alpha, K_{min}}$.  We follow the approach of using regular link functions $\Phi$ as in \cite{NVW18}. 
	
\begin{condition}\label{Cond:LinkFunction1}%-----------------------------------------------------------------

	For given $K_{min}>0$, let $\Phi:\R\to(K_{min},\infty)$ be a smooth, strictly increasing bijective function such that $\Phi(0)=1$, $\Phi'(t)>0, \ t\in\R$, and assume that all derivatives of $\Phi$ are bounded on $\R$.
		
\end{condition}%-------------------------------------------------------------------------------------------------------

	For some of the results to follow it will prove convenient to slightly strengthen the previous condition.

\begin{condition}\label{Cond:LinkFunction2}%-----------------------------------------------------------------

	Let $\Phi$ be as in Condition \ref{Cond:LinkFunction1}, and assume furthermore that $\Phi'$ is nondecreasing and that $\liminf_{t \to -\infty}\Phi'(t)t^a>0$ for some $a>0$.

\end{condition}%-------------------------------------------------------------------------------------------------------

	For $a=2$, an example of such a link function is given in Example \ref{Ex:LinkFunction}  below. Note however that the choice of $\Phi = \exp$ is not permitted in either condition.

	Given any link function $\Phi$ satisfying Condition \ref{Cond:LinkFunction1}, one can show (cf.~\cite{NVW18}, Section 3.1)  that the set $\Fcal_{\alpha,K_{min}}$ in \eqref{Eq:ParamSp} can be realised as the family of composition maps 
$$
	\Fcal_{\alpha,K_{min}}=\{\Phi\circ F : F\in H^\alpha_c(\Ocal)\}, ~~ \alpha \in \mathbb N.
$$
We then regard the solution map associated to \eqref{Eq:DivFormPDE} as one defined on $H^\alpha_c$ via
\begin{equation} 
\label{Eq:ScriptG}
	\Gscr: H^\alpha_c(\Ocal) \to L^2(\Ocal), 
	\quad
		 F\mapsto \Gscr (F):= G(\Phi\circ F),
\end{equation}
where $G(\Phi\circ F)$ is the solution to \eqref{Eq:DivFormPDE} now with $f=\Phi\circ F \in \Fcal_{\alpha, K_{min}}$. In the results to follow, we will implicitly assume a link function $\Phi$ to be given and fixed, and understand the re-parametrised solution map $\Gscr$ as being defined as in \eqref{Eq:ScriptG} for such choice of $\Phi$.

%
%
%

%---------------------------------------------------------------------------------------------------------------------------
%---------------------------------------------------------------------------------------------------------------------------
\subsubsection{Measurement model}
\label{Subsec:MeasModel}
%---------------------------------------------------------------------------------------------------------------------------
%---------------------------------------------------------------------------------------------------------------------------

	Define the uniform distribution on $\Ocal$ by $\mu=dx/\textnormal{vol}(\Ocal)$, where $dx$ is the Lebesgue measure and $\textnormal{vol}(\Ocal)=\int_\Ocal dx$, and consider random design variables 
\begin{equation}
\label{Eq:DiscrRandDesign}
	(X_i)_{i=1}^N \iid \mu, \quad N\in\N.
\end{equation}
For unknown $f\in \Fcal_{\alpha,K_{min}}$, we model the statistical errors under which we observe the corresponding measurements $\{G(f)(X_i)\}_{i=1}^N$ by i.i.d.~Gaussian random variables $W_i\sim N(0,1)$, all independent of the $X_i$'s. Using the re-parameterisation $f= \Phi\circ F$ via a given link function from the previous subsection, the observation scheme is then
\begin{equation}
\label{Eq:DiscrObs}
Y_i=\mathscr G(F)(X_i)+\sigma W_i, \quad i=1,\dots,N,
\end{equation} 
where $\sigma>0$ is the noise amplitude. We will often use the shorthand notation $Y^{(N)}=(Y_i)_{i=1}^N$, with analogous definitions for $X^{(N)}$ and $W^{(N)}$. The random vectors $(Y_i,X_i)$ on $\R\times\Ocal$ are then i.i.d with laws denoted as $P^i_{F}$. Writing $dy$ for the Lebesgue measure on $\R$, it follows that $P^i_{F}$ has Radon-Nikodym density
\begin{equation}
\label{Eq:SingleLikelihood}
	p_F(y,x)
	:=
		\frac{d P^i_{F}}{dy\times d\mu}(y,x)
	=
		\frac{1}{\sqrt{2\pi\sigma^2}}e^{-[y- \mathscr G(F)(x)]^2/(2\sigma^2)},		
	\quad 
		y\in\R,\ x\in\Ocal.
\end{equation}
We will write $P^{N}_{F}=\otimes_{i=1}^N P^i_{F}$ for the joint law of $(Y^{(N)},X^{(N)})$ on $\R^N\times \Ocal^N$, with  $E^i_{F}$, $E^{N}_{F}$ the expectation operators corresponding to the laws $P_F^i,$ $P_F^N$ respectively. In the sequel we sometimes use the notation $P^N_{f}$ instead of $P^N_{F}$ when convenient. 

%
%
%

%---------------------------------------------------------------------------------------------------------------------------
%---------------------------------------------------------------------------------------------------------------------------
\subsubsection{The Bayesian approach}
\label{Subsec:BayesAppr}
%---------------------------------------------------------------------------------------------------------------------------
%---------------------------------------------------------------------------------------------------------------------------

In the Bayesian approach one models the parameter $F \in H^\alpha_c(\Ocal)$ by a Borel probability measure $\Pi$ supported in the Banach space $C(\Ocal)$. Since the map $(F,(y,x))\mapsto p_{F}(y,x)$ can be shown to be jointly measurable, the posterior distribution $\Pi(\cdot|Y^{(N)},X^{(N)})$ of $F|(Y^{(N)},X^{(N)})$ arising from data in model (\ref{Eq:DiscrObs}) equals, by Bayes' formula (p.7, \cite{GvdV17}), 
\begin{equation}
\label{Eq:PostDistr}
	\Pi(B|Y^{(N)},X^{(N)})
	=
		\frac{\int_B e^{\ell^{(N)}(F)}d\Pi(F)}
		{\int_{C(\Ocal)}	e^{\ell^{(N)}(F)}d\Pi(F)}
	\quad
		\textnormal{any Borel set $B\subseteq C(\Ocal)$},
\end{equation}
where
\begin{equation}
\label{Eq:JointLogLikelihood}
	\ell^{(N)}(F)=-\frac{1}{2\sigma^2}\sum_{i=1}^N [Y_i-\Gscr (F)(X_i)]^2
\end{equation}
is (up to an additive constant) the joint log-likelihood function.  

%---------------------------------------------------------------------------------------------------------------------------
%---------------------------------------------------------------------------------------------------------------------------
\subsection{Statistical convergence rates}
\label{Sec:Results}
%---------------------------------------------------------------------------------------------------------------------------
%---------------------------------------------------------------------------------------------------------------------------

	In this section we will show that the posterior distribution arising from certain priors concentrates near any sufficiently regular ground truth $F_0$ (or, equivalently, $f_0$), and provide a bound on the rate of this contraction, assuming the observation $(Y^{(N)},X^{(N)})$ to be generated through model \eqref{Eq:DiscrObs} of law $P_{F_0}^N$. We will regard $\sigma>0$ as a fixed and known constant; in practice it may be replaced by the estimated sample variance of the $Y_i$'s.

	The priors we will consider are built around a Gaussian process base prior $\Pi'$, but to deal with the non-linearity of the inverse problem, some additional regularisation will be required. We first show how this can be done by a $N$-dependent `rescaling' step as suggested in \cite{MNP19b}. We then further show that a randomised truncation of a Karhunen-Loeve-type series expansion of the base prior also leads to a consistent, `fully Bayesian' solution of this inverse problem.
%
%
%

%---------------------------------------------------------------------------------------------------------------------------
%---------------------------------------------------------------------------------------------------------------------------
\subsubsection{Results with re-scaled Gaussian priors}
\label{Subsec:GaussianPriors}
%---------------------------------------------------------------------------------------------------------------------------
%---------------------------------------------------------------------------------------------------------------------------

We will freely use terminology from the basic theory of Gaussian processes and measures, see, e.g., \cite{GN16}, Chapter 2 for details.

\begin{condition}\label{Cond:BasePrior1}%--------------------------------------------------------------------

	Let $\alpha>1+d/2$, $\beta\ge1$, and let $\Hcal$ be a Hilbert space continuously imbedded into $H^\alpha_c(\Ocal)$. Let $\Pi'$ be a centred Gaussian Borel probability measure on the Banach space $C(\Ocal)$ that is supported on a separable measurable linear subspace of $C^\beta(\Ocal)$,
and assume that the reproducing-kernel Hilbert space (RKHS) of $\Pi'$ equals $\Hcal$.

\end{condition}%-------------------------------------------------------------------------------------------------------

As a basic example of a Gaussian base prior $\Pi'$ satisfying Condition \ref{Cond:BasePrior1}, consider a Whittle-Mat\'ern process $M=\{M(x),\ x\in\Ocal\}$ indexed by $\Ocal$ and of regularity $\alpha$ (cf.~Example \ref{Ex:CuttedMaternProcess}  below for full details). We will assume that it is known that $F_0\in H^\alpha(\Ocal)$ is supported inside a given compact subset $K$ of the domain $\Ocal$, and fix any smooth cut-off function $\chi\in C^\infty_c(\Ocal)$ such that $\chi=1$ on $K$. Then, $\Pi'=\Lcal(\chi M)$ is supported on the separable linear subspace $C^{\beta'}(\Ocal)$ of $C^\beta(\Ocal)$ for any $\beta<\beta'<\alpha-d/2$, and its RKHS $\Hcal=\{\chi F, F\in H^\alpha(\Ocal)\}$ is continuously imbedded into $H^\alpha_c(\Ocal)$ (and contains $H^\alpha_K(\Ocal)$). The condition $F_0 \in \mathcal H$ that is employed in the following theorems then amounts to the standard assumption that $F_0 \in H^\alpha(\mathcal O)$ be supported in a strict subset $K$ of $\Ocal$.

\smallskip

	To proceed, if $\Pi'$ is as above and $F'\sim\Pi'$, we consider the `re-scaled' prior
\begin{equation}
\label{Eq:Prior1}
	\Pi_N = \mathcal L(F_N), 
	\quad  F_N = \frac{1}{N^{d/(4\alpha+4+2d)}}F',
\end{equation}
 Then $\Pi_N$ again defines a centred Gaussian prior on $C(\Ocal)$, and a basic calculation (e.g., Exercise 2.6.5 in \cite{GN16}) shows that its RKHS $\Hcal_N$ is still given by $\mathcal H$ but now with norm
\begin{equation}
\label{Eq:RKHS1}
	\|F\|_{\Hcal_N} = N^{d/(4\alpha+4+2d)} \|F\|_\Hcal \quad
	 \forall F\in  \Hcal .
\end{equation}

Our first result shows that the posterior contracts towards $F_0$ in `prediction'-risk at rate $N^{-(\alpha+1)/(2\alpha+2+d)}$ and that, moreover, the posterior draws possess a bound on their $C^\beta$-norm with overwhelming frequentist probability.
	 	
\begin{theorem}\label{Theo:FwdRates1}%----------------------------------------------------------------------

	For fixed integer $\alpha>\beta+d/2$, $\beta\ge1$, consider the Gaussian prior $\Pi_N$ in \eqref{Eq:Prior1} with base prior $F'\sim\Pi'$ satisfying Condition \ref{Cond:BasePrior1} for RKHS $\mathcal H$.  Let $\Pi_N(\cdot|Y^{(N)},X^{(N)})$ be the resulting posterior distribution arising from observations $(Y^{(N)},$ $X^{(N)})$ in \eqref{Eq:DiscrObs}, set $\delta_N=N^{-(\alpha+1)/(2\alpha+2+d)}$, and assume $F_0\in\Hcal$.
	
	Then for any $D>0$ there exists $L>0$ large enough (depending on $\sigma,F_0,D, \alpha, \beta$, as well as on $\Ocal,d,g$) such that, as $N \to \infty$,
\begin{equation}
\label{Eq:FwdRates1}
	\Pi_N(F:\|\Gscr (F)-\Gscr (F_0)\|_{L^2}>L\delta_N|Y^{(N)},X^{(N)}) =O_{P^{N}_{F_0}}(e^{-DN\delta_N^2}),
	\end{equation}
and for sufficiently large $M>0$ (depending on $\sigma, D,\alpha, \beta$)
\begin{equation}
\label{Eq:Regularisation10}
	\Pi_N(F: \|F\|_{C^\beta}>M|Y^{(N)},X^{(N)})=O_{P^{N}_{F_0}}(e^{-DN\delta_N^2}).
\end{equation}
\end{theorem}%------------------------------------------------------------------------------------------------------- 

\smallskip

Following ideas in \cite{MNP19b}, we can combine (\ref{Eq:FwdRates1}) with the regularisation property \eqref{Eq:Regularisation10} and a suitable stability estimate for $G^{-1}$ to show that the posterior contracts about $f_0$ also in $L^2$-risk. We shall employ the stability estimate proved in \cite[Lemma 24]{NVW18}  which requires the source function $g$ in the base PDE (\ref{Eq:DivFormPDE}) to be strictly positive, a natural condition ensuring injectivity of the map $f \mapsto  G(f)$, see  \cite{R81}. Denote the push-forward posterior on the conductivities $f$ by
\begin{equation}
\label{Eq:TildePi}
	\tilde\Pi_N(\cdot|Y^{(N)},X^{(N)}):=\Lcal(f), 
	\quad
	 f=\Phi\circ F : F\sim\Pi_N(\cdot|Y^{(N)},X^{(N)}).
\end{equation}

\begin{theorem}\label{Theo:L2Rates1}%-----------------------------------------------------------------------

	Let $\Pi_N(\cdot|Y^{(N)},X^{(N)})$, $\delta_N$ and $F_0$ be as in Theorem \ref{Theo:FwdRates1} for integer $\beta>1$. Let  $f_0=\Phi\circ F_0$ and assume in addition that $\inf_{x \in \Ocal}g(x) \ge g_{min}>0$. Then for any $D>0$ there exists $L>0$ large enough (depending on $\sigma, f_0, D, \alpha, \beta, \Ocal$, $g_{min}, d$) such that, as $N \to \infty$,
$$
	\tilde\Pi_N(f:\|f-f_0\|_{L^2}>LN^{-\lambda}|Y^{(N)},X^{(N)})
	=
		O_{P^{N}_{f_0}}(e^{-DN\delta_N^2}),
	 \quad
	\lambda
	=
		\frac{(\alpha+1)(\beta-1)}{(2\alpha+2+d)(\beta+1)}.
$$

\end{theorem}%-------------------------------------------------------------------------------------------------------

	We note that as the smoothness $\alpha$ of $f_0$ increases, we can employ priors of higher regularity $\alpha, \beta$. In particular, if $F_0\in C^\infty=\cap_{\alpha>0} H^\alpha$, we can let the above rate $N^{-\lambda}$ be as closed as desired to the `parametric' rate $N^{-1/2}$.

	We conclude this section showing that the posterior mean $ E^\Pi[F|Y^{(N)},X^{(N)}]$ of $\Pi_N(\cdot|Y^{(N)},X^{(N)})$ converges to $F_0$ at the rate $N^{-\lambda}$ from Theorem \ref{Theo:L2Rates1}. We formulate this result at the level of the vector space valued parameter $F$ (instead of for conductivities $f$), as the most commonly used MCMC algorithms (such as pCN, see Subsection \ref{Rem:Computation1}) target the posterior distribution of $F$. 		
\begin{theorem}\label{Theo:PostMeanConv1}%--------------------------------------------------------------

	Under the hypotheses of Theorem \ref{Theo:L2Rates1}, let $\bar F_N=E^\Pi[F|Y^{(N)},X^{(N)}]$ be the (Bochner-) mean of $\Pi_N(\cdot|Y^{(N)},X^{(N)})$. Then, as $N \to \infty$,
\begin{equation}
\label{Eq:PostMeanConv1}
	P^N_{F_0}\big(\|\bar F_N-F_0\|_{L^2}>N^{-\lambda}\big)\to0.
\end{equation}

\end{theorem}

The same result holds for the implied conductivities, that is, for $\|\Phi \circ \bar F_N - f_0\|_{L^2}$ replacing $\|\bar F_N-F_0\|_{L^2}$, since composition with $\Phi$ is Lipschitz.

%--------------------------------------------------------------------------------------------------------

%---------------------------------------------------------------------------------------------------------

\subsubsection{Extension to high-dimensional Gaussian sieve priors}\label{sieve}

It is often convenient, for instance for computational reasons as discussed in Subsection \ref{Rem:Computation1}, to employ  `sieve'-priors that are concentrated on a finite-dimensional approximation of the parameter space supporting the prior.  For example a truncated Karhunen-Loeve-type series expansion (or some other discretisation) of the Gaussian base prior $\Pi'$ is frequently used \cite{DS11, HSV14}. The theorems of the previous subsection remain valid if the approximation spaces are appropriately chosen. 

Let us illustrate this by considering a Gaussian series prior based on an orthonormal basis $\{\Psi_{\ell r}, \ \ell\ge-1,\ r\in\Z^d\}$ of $L^2(\R^d)$ composed of sufficiently regular, compactly supported Daubechies wavelets (see Chapter 4 in \cite{GN16} for details).  We assume that $F_0 \in H^\alpha_K(\Ocal)$ for some $K \subset \Ocal$, and denote by $\Rcal_\ell$ the set of indices $r$ for which the support of $\Psi_{\ell r}$ intersects $K$. Fix any compact $K'\subset \Ocal$ such that $K\subsetneq K'$, and a cut-off function $\chi\in C^\infty_c(\Ocal)$ such that $\chi=1$ on $K'$. For any real $\alpha>1+d/2$, consider the prior $\Pi'_J$ arising as the law of the Gaussian random sum
\begin{equation}
\label{Eq:Prior02}
	\Pi'_J=\Lcal(\chi F), 
	\quad
		 F=\sum_{\ell\le J,r\in\Rcal_\ell}2^{-\ell\alpha}F_{\ell r}\Psi_{\ell r},\ F_{\ell r}\iid N(0,1),
\end{equation}
where $J=J_N \to \infty$ is a (deterministic) truncation point to be chosen. Then $\Pi'_J$ defines a centred Gaussian prior that is supported on the finite-dimensional space 
\begin{equation}
\label{Eq:SubspaceHj}
	\Hcal_J=\textnormal{span}\{\chi \Psi_{\ell r}, \ \ell\le J, \ r\in\Rcal_\ell\}\subset C(\Ocal).
\end{equation}
 
 \begin{proposition}
 Consider a prior $\Pi_N$ as in (\ref{Eq:Prior1}) where now $F' \sim \Pi'_J$ and $J=J_N \in \mathbb N$ is such that $2^{J} \simeq N^{1/(2\alpha +2 +d)}$. Let $\Pi_N(\cdot|Y^{(N)}, X^{(N)})$ be the resulting posterior distribution arising from observations $(Y^{(N)},X^{(N)})$ in \eqref{Eq:DiscrObs}, and assume $F_0 \in H^\alpha_K(\Ocal)$. Then the conclusions of Theorems \ref{Theo:FwdRates1}-\ref{Theo:PostMeanConv1} remain valid (under the respective hypotheses on $\alpha, \beta, g$).
 \end{proposition}
 
A similar result could be proved for more general Gaussian priors (not of wavelet type), but we refrain from giving these extensions here.

%---------------------------------------------------------------------------------------------------------------------------
%---------------------------------------------------------------------------------------------------------------------------
\subsubsection{Randomly truncated Gaussian series priors}
\label{Subsec:RandSeriesPriors}
%--------------------------------------------------------------------------------------------------------------------------

	In this section we show that instead of rescaling Gaussian base priors $\Pi', \Pi'_J$ in a $N-$dependent way to attain extra regularisation, one may also randomise the dimensionality parameter $J$ in (\ref{Eq:Prior02}) by a hyper-prior with suitable tail behaviour. While this is computationally somewhat more expensive (by necessitating a hierarchical sampling method, see Subsection \ref{Rem:Computation1}), it gives a possibly more principled approach to (`fully') Bayesian regularisation in our inverse problem. The theorem below will show that such a procedure is consistent in the frequentist sense, at least for smooth enough $F_0$.

	For the wavelet basis and cut-off function $\chi$ introduced before (\ref{Eq:Prior02}), we consider again a random (conditionally Gaussian) sum
\begin{equation}
\label{Eq:Prior2}
	\Pi=\Lcal(\chi F), 
	\quad
		 F=\sum_{\ell\le J,r\in\Rcal_\ell}2^{-\ell\alpha}F_{\ell r}\Psi_{\ell r},\ F_{\ell r}\iid N(0,1)
\end{equation}
where now $J$ is a random truncation level, independent of the random coefficients $F_{\ell r}$, satisfying the following inequalities
\begin{equation}
\label{Eq:PropertiesOfJ}
	\Pr(J> j)=e^{-2^{jd}\log 2^{jd}}\ \forall j\ge1; 
	\quad 
		\Pr(J=j)\gtrsim  e^{-2^{jd}\log 2^{jd}},\ j\to\infty.
\end{equation}
When $d=1$, a (log-) Poisson random variable satisfies these tail conditions, and for $d>1$ such a random variable $J$ can be easily constructed too -- see Example \ref{poissond} below.

	Our first result in this section shows that the posterior arising from the truncated series prior in \eqref{Eq:Prior2} achieves (up to a log-factor) the same contraction rate in $L^2$-prediction risk as the one obtained in Theorem \ref{Theo:FwdRates1}. Moreover, as is expected in light of the results in \cite{vdVvZ09, R13}, the posterior adapts to the unknown regularity $\alpha_0$ of $F_0$ when it exceeds the base smoothness level $\alpha$.

\begin{theorem}\label{Theo:FwdRates2}%------------------------------------------------------------------

	For any $\alpha>1+d/2$, let $\Pi$ be the random series prior in \eqref{Eq:Prior2}, and let $\Pi(\cdot|Y^{(N)}, X^{(N)})$ be the resulting posterior distribution arising from observations $(Y^{(N)},X^{(N)})$ in \eqref{Eq:DiscrObs}. Then, for each $\alpha_0\ge\alpha$ and any $F_0\in H^{\alpha_0}_K(\Ocal)$, we have that for any $D>0$ there exists $L>0$ large enough (depending on $\sigma,F_0,D,\alpha,\Ocal,d,g$) such that, as $N \to \infty$,
$$
	\Pi(F:\|\Gscr (F)-\Gscr (F_0)\|_{L^2}>L\xi_N|Y^{(N)},X^{(N)})
	=
		O_{P^{N}_{F_0}}(e^{-DN\xi_N^2}),
$$
where $\xi_N= N^{-(\alpha_0+1)/(2\alpha_0+2+d)}\log N$. Moreover, for $\Hcal_J$ the finite-dimensional subspaces in \eqref{Eq:SubspaceHj} and $J_N\in\N$ such that $ 2^{J_N}\simeq N^{1/(2\alpha_0+2+d)}$, we also have that for sufficiently large $M>0$ (depending on $D,\alpha$)
\begin{equation}
\label{Eq:Regularisation2}
	\Pi(F: F\in\Hcal_{J_N}, \ \|F\|_{H^\alpha}\le M 2^{J_N\alpha}N\xi_N^2 |Y^{(N)},X^{(N)})=1-O_{P^{N}_{F_0}}(e^{-DN\xi_N^2}).
\end{equation}

\end{theorem}%--------------------------------------------------------------------------------------------------------

	We can now exploit the previous result along with the finite-dimensional support of the posterior and again the stability estimate from \cite{NVW18}  to obtain the following consistency theorem for $F_0 \in H^{\alpha_0}$ if  $\alpha_0$ is large enough (with a precise bound $\alpha_0 \ge \alpha^*$ given in the proof of Lemma \ref{Lemma:HalphaBound}).
	
\begin{theorem}\label{Theo:L2Rates2}%-----------------------------------------------------------------------

	Let the link function $\Phi$ in the definition \eqref{Eq:ScriptG} of $\Gscr$ satisfy Condition \ref{Cond:LinkFunction2}. Let $\Pi(\cdot|Y^{(N)},X^{(N)})$,  $\xi_N$ be as in Theorem \ref{Theo:FwdRates2}, assume in addition $g \ge g_{min}>0$ on $\Ocal$, and let $\tilde\Pi(\cdot|Y^{(N)},X^{(N)})$ be the posterior distribution of $f$ as in \eqref{Eq:TildePi}. Then for $f_0=\Phi\circ F_0$ with $F_0\in H^{\alpha_0}_K(\Ocal)$ for $\alpha_0$ large enough (depending on $\alpha, d, a$) and for any  $D>0$ there exists $L>0$ large enough (depending on $\sigma,f_0, D, \alpha,\Ocal, g_{min}, d$) such that, as $N \to \infty$,
$$
	\tilde\Pi(f:\|f-f_0\|_{L^2}>LN^{-\rho}|Y^{(N)},X^{(N)}) = O_{P^N_{f_0}}(e^{-DN\xi_N^2}), \quad
	\rho =	\frac{(\alpha_0+1)(\alpha-1)}{(2\alpha_0+2+d)(\alpha+1)}.
$$

\end{theorem}%-------------------------------------------------------------------------------------------------------

	Just as before, for $f_0\in C^\infty$ the above rate can be made as close as desired to $N^{-1/2}$ by choosing $\alpha$ large enough. Moreover, the last contraction theorem also translates into a convergence result for the posterior mean of $F$.

\begin{theorem}\label{Theo:PostMeanConv2}%--------------------------------------------------------------

	Under the hypotheses of Theorem \ref{Theo:L2Rates2}, let $\bar F_N=E^\Pi[F|Y^{(N)},X^{(N)}]$ be the mean of $\Pi(\cdot|Y^{(N)},X^{(N)})$. Then, as $N \to \infty$,
\begin{equation}
\label{Eq:PostMeanConv1}
	P_{F_0}^N \big(\|\bar F_N-F_0\|_{L^2}>N^{-\rho} \big) \to 0.
\end{equation}

\end{theorem}

We note that the proof of the last two theorems crucially takes advantage of the `non-symmetric' and `non-exponential' nature of the stability estimate from \cite{NVW18}, and may not hold in other non-linear inverse problems where such an estimate may not be available (e.g., as in \cite{MNP19b, AN19} or also in the Schr\"odinger equation setting studied in \cite{N17, NVW18}).

\smallskip

Let us conclude this section by noting that hierarchical priors such as the one studied here are usually devised to `adapt to unknown' smoothness $\alpha_0$ of $F_0$, see \cite{vdVvZ09, R13}. Note that while our posterior distribution is adaptive to $\alpha_0$ in the `prediction risk' setting of Theorem \ref{Theo:FwdRates2}, the rate $N^{-\rho}$ obtained in Theorems \ref{Theo:L2Rates2} and \ref{Theo:PostMeanConv2} for the inverse problem \textit{does} depend on the minimal smoothness $\alpha$, and is therefore \textit{not adaptive}. Nevertheless, this hierarchical prior gives an example of a fully Bayesian, consistent solution of our inverse problem.

\subsection{Concluding discussion}

\subsubsection{Posterior computation}\label{Rem:Computation1}\normalfont%-----------------------------------------------------
As mentioned in the introduction, in the context of the elliptic inverse problem considered in the present paper, posterior distributions arising from Gaussian process priors such as those above can be computed by MCMC algorithms, see \cite{CRSW13, CMPS16, BGLFS17}, and computational guarantees can be obtained as well: For Gaussian priors, \cite{HSV14} establish non-asymptotic sampling bounds for the `preconditioned Crank-Nicholson (pCN)' algorithm, which hold even in the absence of log-concavity of the likelihood function, and which imply bounds on the approximation error for the computation of the posterior mean.  The algorithm can be implemented as long as it is possible to evaluate the forward map $F\mapsto \Gscr(F)(x)$ at $x \in \mathcal O$, which in our context can be done by using standard numerical methods to solve the elliptic PDE \eqref{Eq:DivFormPDE}. In practice, these algorithms often employ a finite-dimensional approximation of the parameter space (see Subsection \ref{sieve}).
	
	In order to sample from the posterior distribution arising from the more complex hierarchical prior \eqref{Eq:Prior2}, MCMC methods based on fixed Gaussian priors (such as the pCN algorithm) can be employed within a suitable Gibbs-sampling scheme that exploits the conditionally Gaussian structure of the prior. The algorithm would then alternate, for given $J$, an MCMC step targeting the marginal posterior distribution of $F|(Y^{(N)},X^{(N)},J)$, followed by, given the actual sample of $F$, a second MCMC run with objective the marginal posterior of $J|(Y^{(N)},X^{(N)},F)$. A related approach to hierarchical inversion is empirical Bayesian estimation. In the present setting this would entail first estimating the truncation level $J$ from the data, via an estimator $\hat J=\hat J(Y^{(N)},X^{(N)})$ (e.g., the marginal maximum likelihood estimator), and then performing inference based on the fixed finite-dimensional prior $\Pi_{\hat J}$ (defined as in \eqref{Eq:Prior2} with $J$ replaced by $\hat J$). See \cite{K2015} where this is studied in a diagonal linear inverse problem.

\subsubsection{Open problems: Towards optimal convergence rates}

The convergence rates obtained in this article  demonstrate the frequentist consistency of a Bayesian (Gaussian process) inversion method in the elliptic inverse problem (\ref{Eq0}) with data (\ref{discrete}) in the large sample limit $N \to \infty$. While the rates approach the optimal rate $N^{-1/2}$ for very smooth models ($\alpha \to \infty$), the question of optimality for fixed $\alpha$ remains an interesting avenue for future research.  We note that for the `PDE-constrained regression' problem of recovering $\mathscr G(F_0)$ in `prediction' loss, the rate $\delta_N=N^{-(\alpha+1)/(2\alpha+2+d)}$ obtained in Theorems \ref{Theo:FwdRates1} and \ref{Theo:FwdRates2} can be shown to be minimax optimal (as in \cite[Theorem 10]{NVW18}). But for the recovery rates for $f$ obtained in Theorems \ref{Theo:PostMeanConv1} and \ref{Theo:PostMeanConv2}, no matching lower bounds are currently known. Related to this issue, in \cite{NVW18} faster (but still possibly suboptimal) rates are obtained for the modes of our posterior distributions (MAP estimates, which are not obviously computable in polynomial time), and one may loosely speculate here about computational hardness barriers in our non-linear inverse problem. These issues pose formidable challenges for future research and are beyond the scope of the present paper.

%------------------------------------------------------------------------------------------------------------------------
%------------------------------------------------------------------------------------------------------------------------
\section{Proofs}
\label{Sec:Proofs}
%------------------------------------------------------------------------------------------------------------------------
%------------------------------------------------------------------------------------------------------------------------

We assume without loss of generality that $\vol(\Ocal)=1$. In the proof, we will repeatedly exploit properties of the (re-parametrised) solution map $\Gscr$ defined in \eqref{Eq:ScriptG}, which was studied in detail in \cite{NVW18}. Specifically, in the proof of Theorem 9 in \cite{NVW18} it is shown that, for all $\alpha>1+d/2$ and any $F_1, F_2\in H^\alpha_c(\Ocal)$, 
\begin{equation}
\label{Eq:LipCondG}
	\|\Gscr (F_1)-\Gscr (F_2)\|_{L^2(\Ocal)}
	\lesssim 
		(1+\|F_1\|_{C^1(\Ocal)}^4\vee\|F_2\|_{C^1(\Ocal)}^4)\|F_1-F_2\|_{(H^1(\Ocal))^*},
\end{equation}
where we denote by $X^*$ the topological dual Banach space of a normed linear space $X$. Secondly, we have (Lemma 20 in \cite{NVW18}) for some constant $c>0$ (only depending on $d,\ \Ocal$ and $K_{min})$, 
\begin{equation}
\label{Eq:UnifBoundG}
	\sup_{F\in H^\alpha_c}\|\Gscr (F)\|_{\infty}\le c\|g\|_{\infty}<\infty.
\end{equation}
Therefore the inverse problem \eqref{Eq:DiscrObs} falls in the general framework considered in Appendix \ref{App:GenInvProbl} below (with $\beta=\kappa=1$, $\gamma=4$ in \eqref{Eq:GenLipCondG} and $S=c\|g\|_\infty$ in \eqref{Eq:GenUnifBoundG}) ; in particular Theorems \ref{Theo:FwdRates1} and \ref{Theo:FwdRates2} then follow as particular cases of the general contraction rate results derived in Theorem \ref{Theo:GenFwdRates1} and Theorem \ref{Theo:GenFwdRates2}, respectively. It thus remains to derive Theorems \ref{Theo:L2Rates1} and \ref{Theo:PostMeanConv1} from Theorem \ref{Theo:FwdRates1}, and Theorems \ref{Theo:L2Rates2} and \ref{Theo:PostMeanConv2} from Theorem \ref{Theo:FwdRates2}, respectively.

To do so we recall here another key result from \cite{NVW18}, namely their stability estimate Lemma 24: For $\alpha>2+d/2$, if $G(f)$ denotes the solution of the PDE (\ref{Eq:DivFormPDE}) with $g$ satisfying $\inf_{x\in\Ocal}g(x)\ge g_{min}>0$, then for fixed $f_0\in \Fcal_{\alpha,K_{min}}$ and all $f\in\Fcal_{\alpha,K_{min}}$
\begin{equation}
\label{Eq:StabEstim}
	\|f-f_0\|_{L^2(\Ocal)}\lesssim \|f\|_{C^1(\Ocal)}\|G(f)-G(f_0)\|_{H^2(\Ocal)},
\end{equation}
with multiplicative constant independent of $f$.
%
%
%

%------------------------------------------------------------------------------------------------------------------------
%------------------------------------------------------------------------------------------------------------------------
\subsection{Proofs for Section \ref{Subsec:GaussianPriors}}
\label{Subsec:ProofsGaussPriors}
%------------------------------------------------------------------------------------------------------------------------
%------------------------------------------------------------------------------------------------------------------------

\paragraph{Proof of  Theorem \ref{Theo:L2Rates1}.}
%------------------------------------------------------------------------------------------------------------------------

	The conclusions of Theorem \ref{Theo:FwdRates1} can readily be translated for the push-forward posterior $\tilde\Pi_N(\cdot|Y^{(N)},X^{(N)})$ from \eqref{Eq:TildePi}. In particular, \eqref{Eq:FwdRates1} implies, for $f_0=\Phi\circ F_0$, as $N\to\infty$,
\begin{equation}
\label{Eq:L2FwdContrRate1}
	\tilde\Pi_N(f : \|G(f)-G(f_0)\|_{L^2}>L \delta_N |Y^{(N)},X^{(N)})
	=
		O_{P^N_{f_0}}(e^{-D N\delta_N^2 });
\end{equation}
and using Lemma 29 in \cite{NVW18} and \eqref{Eq:Regularisation10} we obtain for sufficiently large $M'>0$
 \begin{equation}
 \label{Eq:NormBound1}
	\tilde\Pi_N(f : \|f\|_{C^\beta}>M'|Y^{(N)},X^{(N)})
	\le
		\Pi_N( F : \|F\|_{C^{\beta}}>M|Y^{(N)},X^{(N)})\\
	=
		O_{P^N_{f_0}}(e^{-D N\delta_N^2 }).
\end{equation}

	From the previous bounds we now obtain the following result.
\begin{lemma}\label{Lemma:FwdH2Risk1}%-------------------------------------------------------------------

	For $\Pi_N(\cdot|Y^{(N)},X^{(N)}), \delta_N$ and $F_0$ as in Theorem \ref{Theo:FwdRates1}, let $\tilde\Pi_N(\cdot|Y^{(N)},X^{(N)})$ be the push-forward posterior distribution from \eqref{Eq:TildePi}. Then, for $f_0=\Phi\circ F_0$ and any $D>0$ there exists $L>0$ large enough such that, as $N\to\infty$,
$$
	\tilde\Pi_N(f:\|G(f)-G(f_0)\|_{H^2}>L \delta_N ^{(\beta-1)/(\beta+1)}|Y^{(N)},X^{(N)})
	= O_{P^N_{F_0}}(e^{-D N\delta_N^2 }).
$$

\proof%------------------------------------------------------------------------------------------------------------------

	Using the continuous imbedding of $C^\beta \subset H^\beta, \beta \in \mathbb N,$ and \eqref{Eq:NormBound1}, for some $M'>0$ as $N \to \infty$,
$$
	\tilde\Pi_N(f : \|f\|_{H^\beta}>M'|Y^{(N)},X^{(N)})=O_{P^N_{F_0}}(e^{-DN\delta_N^2}).
$$
 Now if $f\in H^\beta$ with $\|f\|_{H^\beta}\le M'$, Lemma 23 in \cite{NVW18} implies $G(f), G(f_0)\in H^{\beta+1}$, with
$$
	\|G(f_0)\|_{H^{\beta+1}}\lesssim 1+\|f_0\|_{H^\beta}^{\beta(\beta+1)}<\infty,	
	\quad
	\|G(f)\|_{H^{\beta+1}}\lesssim 1+\|f\|_{H^\beta}^{\beta(\beta+1)}	
	<M''<\infty;
$$
and by the usual interpolation inequality for Sobolev spaces \cite{Lions1972},
\begin{align*}
	\|G(f)-G(f_0)\|_{H^2}
	&\lesssim 
		\|G(f)-G(f_0)\|^{(\beta-1)/(\beta+1)}_{L^2}\|G(f)-G(f_0)\|
		^{2/(\beta+1)}_{H^{\beta+1}}
		\\
%	&\le
%		\|G(f)-G(f_0)\|^{(\beta_0-1)/(\beta_0+1)}_{L^2}(\|G(f)\|_{H^{\beta_0+1}}
%		+\|G(f_0)\|_{H^{\beta_0+1}})^{2/(\beta_0+1)}\\
	&\lesssim
		\|G(f)-G(f_0)\|^{(\beta-1)/(\beta+1)}_{L^2}.
\end{align*}
Thus, by what precedes and \eqref{Eq:L2FwdContrRate1}, for sufficiently large $L>0$
\begin{align*}
	&\tilde\Pi_N(f:\|G(f)-G(f_0)\|_{H^2} 
	> L \delta_N ^{(\beta-1)/(\beta+1)}|Y^{(N)},X^{(N)}) \\
	&\ \le 
		\tilde\Pi_N(f:\|G(f)-G(f_0)\|_{L^2}> L' \delta_N |Y^{(N)},X^{(N)})
		+\tilde\Pi_N(f:\|f\|_{H^{\beta}}> M''|Y^{(N)},X^{(N)})\\
	&\ =
		O_{P^N_{F_0}}(e^{-D N\delta_N^2 }),
\end{align*}
as $N \to \infty$.

\endproof%--------------------------------------------------------------------------------------------------------------

\end{lemma}%----------------------------------------------------------------------------------------------------------

	To prove Theorem \ref{Theo:L2Rates1} we use \eqref{Eq:StabEstim}, \eqref{Eq:NormBound1} and Lemma \ref{Lemma:FwdH2Risk1} to the effect that for any $D>0$ we can find $L, M>0$ large enough such that, as $N \to \infty$,
\begin{align*}
	&\tilde\Pi_N(f:\|f-f_0\|_{L^2}> L \delta_N ^\frac{\beta-1}{\beta+1}|Y^{(N)},X^{(N)})\\
	&\  \le 
		\tilde \Pi_N(f:\|G(f)-G(f_0)\|_{H^2}> L' \delta_N ^\frac{\beta-1}{\beta+1}
		|Y^{(N)},X^{(N)})+\tilde\Pi_N(f:\|f\|_{C^\beta}> M|Y^{(N)},X^{(N)})\\
	&\ = O_{P^N_{F_0}}(e^{-D N\delta_N^2 }).
\end{align*}

%------------------------------------------------------------------------------------------------------------------------
\paragraph{Proof of  Theorem \ref{Theo:PostMeanConv1}.}
%------------------------------------------------------------------------------------------------------------------------

The proof largely follows ideas of \cite{MNP19b} but requires a slightly more involved, iterative uniform integrability argument to also control the probability of events $\{F:\|F\|_{C^\beta} >M\}$ on whose complements we can subsequently exploit regularity properties of the inverse link function $\Phi^{-1}$.

Using Jensen's inequality, it is enough to show, as $N\to\infty$,
\begin{align*}
	P_{F_0}^N\Big(E^\Pi[\|F-F_0\|^2_{L^2}|Y^{(N)},X^{(N)}]>N^{-\lambda}\Big)\to0.
\end{align*}
For $M>0$ sufficiently large to be chosen, we decompose
\begin{align}
\label{Eq:TwoThings}
	 E^\Pi[\|F-F_0\|_{L^2}|Y^{(N)},X^{(N)}]
	&=
		 E^\Pi[\|F-F_0\|_{L^2}1_{\|F\|_{C^\beta}\le M}
		|Y^{(N)},X^{(N)}]\nonumber\\
	&\quad+
		 E^\Pi[\|F-F_0\|_{L^2}1_{\|F\|_{C^\beta}> M}
		|Y^{(N)},X^{(N)}].
\end{align}
Using the Cauchy-Schwarz inequality we can upper bound the expectation in the second summand by
\begin{align*}
&	\sqrt{E^\Pi[\|F-F_0\|^2_{L^2}|Y^{(N)},X^{(N)}]}
		\sqrt{\Pi_N(F:\|F\|_{C^\beta}> M|Y^{(N)},X^{(N)})}.
\end{align*}
In view of \eqref{Eq:Regularisation10}, for all $D>0$ we can choose $M>0$ large enough to obtain
\begin{align*}
	P_{F_0}^N&\Big(E^\Pi[\|F-F_0\|^2_{L^2}|Y^{(N)},X^{(N)}]
		\Pi_N(F:\|F\|_{C^\beta}> M|Y^{(N)},X^{(N)})>N^{-2\lambda}\Big)\\
	&\le
		P_{F_0}^N\Big(E^\Pi[\|F-F_0\|^2_{L^2}|Y^{(N)},X^{(N)}]
		e^{-DN\delta_N^2}>N^{-2\lambda}\Big)+o(1).
\end{align*}
To bound the probability in the last line, let $\Bcal_N$ be the sets defined in \eqref{Eq:Bn} below, note that Lemma \ref{Lemma:SmallBall1} and Lemma \ref{Lemma:HellingerAndL2SmallBalls} below jointly imply that $\Pi_N(\Bcal_N)\ge ae^{-A N\delta_N^2 }$ for some $a,A>0$. Also, let $\nu(\cdot)=\Pi_N(\cdot\cap\Bcal_N)/\Pi_N(\Bcal_N)$, and let $\Ccal_N$ be the event from \eqref{Eq:Cn}, for which Lemma 7.3.2 in \cite{GN16} implies that $P^N_{F_0}(\Ccal_N)\to1$ as $N\to\infty$. Then
\begin{align*}
	&P_{F_0}^N\Big(E^\Pi[\|F-F_0\|^2_{L^2}|Y^{(N)},X^{(N)}]
	e^{-DN\delta_N^2}>N^{-2\lambda}\Big)\\
	&\  \le
		P_{F_0}^N\Bigg(\frac{\int_{C(\Ocal)}\|F-F_0\|^2_{L^2}\prod_{i=1}^Np_F/p_{F_0}(Y_i,X_i)d\Pi_N(F)}
		{\Pi(\Bcal_N)\int_{\Bcal_N}\prod_{i=1}^Np_F/p_{F_0}(Y_i,X_i)d\nu(F)}
		e^{-DN\delta_N^2}>N^{-2\lambda},\Ccal_N\Bigg)+o(1)\\
	&\  \le
		P_{F_0}^N\Big(\int_{C(\Ocal)}\|F-F_0\|^2_{L^2}\prod_{i=1}^Np_F/p_{F_0}(Y_i,X_i)d\Pi_N(F)
		>N^{-2\lambda}ae^{(D-A-2)N\delta_N^2} \Big)+o(1)
\end{align*}
which is upper bounded, using Markov's inequality and Fubini's theorem, by
\begin{align*}
	&\frac{1}{a}e^{-(D-A-2)N\delta_N^2}N^{2\lambda}
		\int_{C(\Ocal)}\|F-F_0\|^2_{L^2}E_{F_0}^N\Bigg(\prod_{i=1}^N\frac{p_F}{p_{F_0}}(Y_i,X_i)\Bigg)d\Pi_N(F).
\end{align*}
Taking $D>A+2$ (and $M$ large enough in \eqref{Eq:TwoThings}), using the fact that $E_{F_0}^N\big(\prod_{i=1}^N$ $p_F/p_{F_0}(Y_i,X_i)\big)=1$, and that $E^{\Pi_N}\|F\|_{L^2}<\infty$ (by Fernique's theorem, e.g., \cite[Exercise 2.1.5]{GN16}), we then conclude
\begin{equation}
\label{idk1}
	P_{F_0}^N\Big(E^\Pi[\|F-F_0\|^2_{L^2}1_{\|F\|_{C^\beta}>M}|Y^{(N)},X^{(N)}]>N^{-\lambda}\Big)\to0, \quad N\to\infty.
\end{equation}

	To handle the first term in \eqref{Eq:TwoThings}, let $f=\Phi\circ F$ and $f_0=\Phi\circ F_0$. Then for all $x\in\Ocal$, by the mean value and inverse function theorems,
\begin{align*}
	|F(x)-F_0(x)|
	=
	|\Phi^{-1}\circ f(x)-\Phi^{-1}\circ f_0(x)| 
%	&=
%		|(\Phi^{-1})'(\eta)||f(x)-f_0(x)|\\
	&=
		\frac{1}{|\Phi'(\Phi^{-1}(\eta))|}|f(x)-f_0(x)|
\end{align*}
for some $\eta$ lying between $f(x)$ and $f_0(x)$. If $\|F\|_{C^\beta}\le M$ then, as $\Phi$ is strictly increasing, necessarily
$
	f(x)=\Phi(F(x))
	\in [\Phi(-M), \Phi(M)]
$
for all $x\in\Ocal$.
Similarly, the range of $f_0$ is contained in the compact interval $[\Phi(-M), \Phi(M)]$ for $M \ge \|F_0\|_\infty$, so that
\begin{align*}
	|\Phi^{-1}\circ f(x)-\Phi^{-1}\circ f_0(x)| 
	&\le
		\frac{1}{\min_{z \in [-M,M]}\Phi'(z)} |f(x)-f_0(x)| \lesssim |f(x)-f_0(x)|
\end{align*}
for a multiplicative constant not depending on $x\in\Ocal$. It follows 
\begin{align*}
	\|F-F_0\|_{L^2} 1_{\|F\|_{C^\beta}\le M}
	&\lesssim
		\|f-f_0\|_{L^2}1_{\|F\|_{C^\beta}\le M},
\end{align*}
and 
$$
	E^\Pi[\|F-F_0\|_{L^2}1_{\|F\|_{C^\beta}\le M}
		|Y^{(N)},X^{(N)}]
	\lesssim
		E^{\tilde\Pi}[\|f-f_0\|_{L^2}
		|Y^{(N)},X^{(N)}].
$$
Noting that for each $L>0$ the last expectation is upper bounded by
\begin{align*}
	L&N^{-\lambda}
		+ E^{\tilde\Pi}\Big[\|f-f_0\|_{L^2}
		1_{\|f-f_0\|_{L^2}>LN^{-\lambda}}|Y^{(N)},X^{(N)}]\\
		&\le
			LN^{-\lambda}+
			\sqrt{E^{\tilde\Pi}[\|f-f_0\|^2_{L^2}|Y^{(N)},X^{(N)}]}
			\sqrt{\tilde\Pi_N(f:\|f-f_0\|_{L^2}>LN^{-\lambda}|Y^{(N)},X^{(N)})},
\end{align*}
we can repeat the above argument, with the event $\{F:\|F\|_{C^\beta} > M\}$ replaced by the event $\{f:\|f-f_0\|_{L^2}>LN^{-\lambda}\}$, to deduce from Theorem \ref{Theo:L2Rates1} that for $D>A+2$ there exists $L>0$ large enough such that
\begin{align*}
	P_{F_0}^N&\Big(E^{\tilde\Pi}[\|f-f_0\|^2_{L^2}|
	Y^{(N)},X^{(N)}]
	\tilde\Pi_N(f:\|f-f_0\|^2_{L^2}>LN^{-\lambda}|Y^{(N)},X^{(N)})
	>N^{-\lambda}\Big)\\
	&\lesssim 
	e^{-(D-A-2)N\delta_N^2}N^{2\lambda}
\end{align*}
which combined with \eqref{idk1} and the definition of $\delta_N$ concludes the proof.

%
%
%

%--------------------------------------------------------------------------------------------------------------------------------------------------------------
%--------------------------------------------------------------------------------------------------------------------------------------------------------------
\subsection{Sieve prior proofs}
%--------------------------------------------------------------------------------------------------------------------------------------------------------------
%--------------------------------------------------------------------------------------------------------------------------------------------------------------

The proof only requires minor modification from the proofs of Section \ref{Subsec:GaussianPriors}. We only discuss here the main points: One first applies the $L^2$-prediction risk Theorem \ref{Theo:GenFwdRates1}  with a sieve prior. In the proof of the small ball Lemma \ref{Lemma:SmallBall1} one uses the following observations: the projection $P_{\Hcal_J}(F_0)\in \Hcal_J$ of $F_0 \in H^\alpha_K$ defined in \eqref{Eq:Projections} satisfies by (\ref{Eq:DualNormApprox})
$$\|F_0-P_{\Hcal_J}(F_0)\|_{(H^1(\Ocal))^*}
	\lesssim	2^{-J(\alpha+1)};$$
hence choosing $J$ such that $2^J\simeq N^{1/(2\alpha+2+d)}$, and noting also that $\|P_{\Hcal_J}(F_0)\|_{C^1}\le \|F_0\|_{C^1}<\infty$ for all $J$ by standard properties of wavelet bases, it follows from \eqref{Eq:LipCondG} that
$$
 	\|\Gscr (F_0)-\Gscr (P_{\Hcal_J}(F_0))\|_{L^2}
	\lesssim
		 \|F_0-P_{\Hcal_J}(F_0)\|_{(H^1)^*}
	\lesssim
		N^{-(\alpha+1)/(2\alpha+2+d)}
	=
		 \delta_N.
$$
Therefore, by the triangle inequality,
$$
	\Pi_N(F:\|\Gscr (F)-\Gscr (F_0)\|_{L^2}\ge \delta_N/q)\ge \Pi_N(F:\|\Gscr (F)-\Gscr (P_{\Hcal_N}(F_0))\|_{L^2}\ge q'\delta_N).
$$
The rest of the proof of Lemma \ref{Lemma:SmallBall1} then carries over (with $P_{\Hcal_J}(F_0)$ replacing $F_0$), upon noting that (\ref{Eq:Norms}) and a Sobolev imbedding imply
$$
	\sup_{J\in\N} E^{\Pi'_J}\|F\|^2_{C^1}<\infty,~\text{ as well as }~\|F\|_{H^\alpha}\le c\|F\|_{\Hcal_J} \text{ for all } F\in\Hcal_J
$$ for some constant $c>0$ independent of $J$. Moreover, the last two properties are sufficient to prove an analogue of Lemma \ref{Lemma:ApproxSets1} as well, so that Theorem \ref{Theo:GenFwdRates1} indeed applies to the sieve prior. The proof from here onwards is identical to the ones of Theorems \ref{Theo:FwdRates1}-\ref{Theo:PostMeanConv1} for the unsieved case, using also that what precedes implies that $\sup_J E^{\Pi'_J}\|F\|^2_{L^2}<\infty$, relevant in the proof of convergence of the posterior mean.

%------------------------------------------------------------------------------------------------------------------------
%------------------------------------------------------------------------------------------------------------------------
\subsection{Proofs for Section \ref{Subsec:RandSeriesPriors}}
\label{Subsec:ProofsRandSeriesPrior}
%------------------------------------------------------------------------------------------------------------------------
%------------------------------------------------------------------------------------------------------------------------

Inspection of the proofs for rescaled priors implies that Theorems \ref{Theo:L2Rates2} and \ref{Theo:PostMeanConv2} can be deduced from Theorem \ref{Theo:FwdRates2} if we can show that posterior draws lie in a $\alpha$-Sobolev ball of fixed radius with sufficiently high frequentist probability. This is the content of the next result.

\begin{lemma}\label{Lemma:HalphaBound}%---------------------------------------------------------------

	Under the hypotheses of Theorem \ref{Theo:L2Rates2}, there exists $\alpha^* >0$ (depending on $\alpha, d$ and $a$) such that for each $F_0\in H^{\alpha_0}_K(\Ocal), \alpha_0>\alpha^*,$ and any $D>0$ we can find $M>0$ large enough such that, as $N \to \infty$,
 $$
	\Pi(F: \|F\|_{H^\alpha}\le M|Y^{(N)},X^{(N)})
	=
		1-O_{P^N_{F_0}}(e^{-DN\xi_N^2}).
$$
 
 \proof%-----------------------------------------------------------------------------------------------------------------
 
	Theorem \ref{Theo:FwdRates2} implies that for all $D>0$ and sufficiently large $L,M>0$, if $J_N\in\N:2^{J_N}\simeq N^{1/(2\alpha_0+2+d)}$ and denoting by
\begin{align*}
	\Acal_N
	&=
		\{F\in \Hcal_{J_N} :  
		\| F\|_{H^\alpha}\le M2^{J_N\alpha}\sqrt{N}\xi_N,
		 \ \|\Gscr(F)-\Gscr (F_0)\|_{L^2}\le L\xi_N\}, 
\end{align*}
then as $N\to \infty$
\begin{equation}
\label{Eq:AvarepsilonDecay}
	 \Pi(\Acal_N|Y^{(N)},X^{(N)})
	 	=1-O_{P^N_{F_0}}(e^{-DN\xi_N^2}).
 \end{equation}

	Next,	 note that if $F\in\Hcal_{J_N}$, then by standard properties of wavelet bases (cf.~\eqref{Eq:EquivOfNorms}),
$
	\| F\|_{H^\alpha}
	\lesssim 
		2^{J_N\alpha}\| F\|_{L^2}
$
 for all $N$ large enough. Thus, for $P_{\Hcal_{J_N}}(F_0)$ the projection of $F_0$ onto $\Hcal_{J_N}$ defined in \eqref{Eq:Projections},
$$
	\|F\|_{H^\alpha}
	\le
		\|F-P_{\Hcal_{J_N}}(F_0)\|_{H^\alpha}+\|P_{\Hcal_{J_N}}(F_0)\|_{H^\alpha}
	\lesssim 
		2^{J_N\alpha}\|F-F_0\|_{L^2}+\|F_0\|_{H^\alpha},
$$
and a Sobolev imbedding further gives $\|F\|_{L^\infty}\le M' 2^{J_N\alpha}\sqrt{N}\xi_N$, for some $M'>0$. Now letting $f=\Phi\circ F$ and $f_0=\Phi\circ F_0$, by similar argument as in the proof of Theorem \ref{Theo:PostMeanConv1} combined with monotonicity of $\Phi'$, we see that for all $N$ large enough
\begin{align*}
	\|F-F_0\|_{L^2}
	\le\frac{1}{\Phi'(-M'2^{J_N\alpha}\sqrt{N}\xi_N)} \|f-f_0\|_{L^2}.
\end{align*}
Then, using the assumption on the left tail of $\Phi$ in Condition \ref{Cond:LinkFunction2}, and the stability estimate \eqref{Eq:StabEstim},
 \begin{align*}
	\|F-F_0\|_{L^2}
	&\lesssim
		(2^{J_N\alpha}\sqrt{N}\xi_N)^a\|f\|_{H^\alpha}\|G(f)-G(f_0)\|_{H^2}.
\end{align*}
Finally, by the interpolation inequality for Sobolev spaces \cite{Lions1972} and Lemma 23 in \cite{NVW18}, 
\begin{align*}
	\|G(f)-G(f_0)\|_{H^2}
	&\lesssim
		\|G(f)-G(f_0)\|_{L^2}^{(\alpha-1)/(\alpha+1)}
		\|G(f)-G(f_0)\|^{2/(\alpha+1)}_{H^{\alpha+1}}\\
	&\lesssim
		\xi_N^{(\alpha-1)/(\alpha+1)} 
		(\|G(f)\|_{H^{\alpha+1}}+\|G(f_0)\|_{H^{\alpha+1}})^{2/(\alpha+1)}\\
	&\lesssim
		\xi_N^{(\alpha-1)/(\alpha+1)}
		 (1+\|f\|^{\alpha^2+\alpha}_{H^{\alpha}})^{2/(\alpha+1)}.
\end{align*}
so that, in conclusion, for each $F\in\Acal_N$ and sufficiently large $N$,
\begin{align*}
	\| F\|_{H^\alpha}
	&\lesssim
		1+
		2^{J_N\alpha}(2^{J_N\alpha}\sqrt{N}\xi_N)^a\|f\|_{H^\alpha}
		\xi_N^{\frac{\alpha-1}{\alpha+1}}
		(1+\|f\|_{H^\alpha}^{\alpha^2+\alpha})^\frac{2}{\alpha+1}.
\end{align*}
The last term is bounded, using Lemma 29 in \cite{NVW18}, by a multiple of 
\begin{align*}
%	\xi_N^\frac{\alpha-1}{\alpha+1}
%	2^{J_N\alpha}
%	(2^{J_N\alpha}\sqrt{N}\xi_N)^a
%	(2^{J_N\alpha}\sqrt{N}\xi_N)^\alpha&
%	(2^{J_N\alpha}\sqrt{N}\xi_N)^{2\alpha^2}\\
%	&=
%		N^{\frac{\alpha}{2\alpha_0+2+d}}
%		\xi_N^\frac{\alpha-1}{\alpha+1}
%		(2^{J_N\alpha}\sqrt{N}\xi_N)^{2\alpha^2+\alpha+a}.
	\xi_N^\frac{\alpha-1}{\alpha+1}
	2^{J_N\alpha}
	(2^{J_N\alpha}\sqrt{N}\xi_N)^{2\alpha^2+2\alpha+a}
	=
	N^{-\frac{(\alpha_0+1)(\alpha-1)}
	{(2\alpha_0+2+d)(\alpha+1)}}
	N^\frac{2\alpha^3+(2+d)\alpha^2+(1+a+d)\alpha+ad/2}{2\alpha_0+2+d}
\end{align*}
the last identity holding up to a log factor. Hence, if
$$
	\alpha_0
	>
		\alpha^*
	:=
		\frac{[2\alpha^3+(2+d)\alpha^2+(1+a+d)\alpha+ad/2](\alpha+1)}{(\alpha-1)}
$$
then we conclude overall that $\|F\|_{H^\alpha}\lesssim 1+o(1)$ as $N\to\infty$ for all $F\in \Acal_N$, proving the claim in view of \eqref{Eq:AvarepsilonDecay}.

\endproof%-------------------------------------------------------------------------------------------------------------

\end{lemma}%---------------------------------------------------------------------------------------------------------

	Replacing $\beta$ by $\alpha$ in the conclusion of Lemma \ref{Lemma:FwdH2Risk1}, the proof of Theorem \ref{Theo:L2Rates2} now proceeds as in the proof of Theorem \ref{Theo:L2Rates1} without further modification. Likewise, Theorem \ref{Theo:PostMeanConv2} can be shown following the same argument as in the proof of Theorem \ref{Theo:PostMeanConv1}, noting that for $\Pi$ the random series prior in \eqref{Eq:Prior2}, it also holds that $E^\Pi\|F\|^2_{L^2}<\infty$.

%
%
%
%
%

%------------------------------------------------------------------------------------------------------------------------
%------------------------------------------------------------------------------------------------------------------------
\appendix
\section{Results for general inverse problems}
\label{App:GenInvProbl}
%------------------------------------------------------------------------------------------------------------------------
%------------------------------------------------------------------------------------------------------------------------

	Let $\Ocal\subset\R^d, \ d\in\N$, be a nonempty open and bounded set with smooth boundary, and assume that $\Dcal$ is a nonempty and bounded measurable subset of $\R^p, \ p\ge1$. Let $\Fcal\subseteq L^2(\Ocal)$ be endowed with the trace Borel $\sigma$-field of $L^2(\Ocal)$, and consider a Borel-measurable ‘forward mapping'
$$
	\Gscr: \Fcal \to L^2(\Dcal), \quad F\mapsto \Gscr (F).
$$

	For $F\in\Fcal$, we are given noisy discrete measurements of $\Gscr(F)$ over a grid of points drawn uniformly at random on $\Dcal$, 
\begin{equation}
\label{Eq:GenDiscrObs}
	Y_i=\Gscr (F)(X_i)+\sigma W_i, \quad i=1,\dots,N,\ X_i\iid\mu, \ W_i\iid N(0,1),
\end{equation} 
for some $\sigma>0$. Above $\mu$ denotes the uniform (probability) distribution on $\Dcal$ and the design variables $(X_i)_{i=1}^N$ are independent of the noise vector $(W_i)_{i=1}^N$. We assume without loss of generality that $\vol(\Dcal)=1$, so that $\mu=dx$, the Lebesgue measure on $\Dcal$.

	We take the noise amplitude $\sigma>0$ in \eqref{Eq:GenDiscrObs} to be fixed and known, and work under the assumption that the forward map $\Gscr$ satisfies the following local Lipschitz condition: for given $\beta,\gamma,\kappa\ge0$, and all $F_1,F_2 \in C^\beta(\Ocal) \cap \mathcal F$,
\begin{equation}
\label{Eq:GenLipCondG}
	\|\Gscr(F_1)-\Gscr(F_2)\|_{L^2(\Dcal)} 
	\lesssim
		 (1+\|F_1\|_{C^\beta(\Ocal)}^\gamma \vee \|F_2\|_{C^\beta(\Ocal)}^\gamma)
		  \|F_1-F_2\|_{(H^\kappa(\Ocal))^*}
\end{equation}
where we recall that $X^*$ denotes the topological dual Banach space of a normed linear space $X$. Additionally, we will require $\Gscr$ to be uniformly bounded on its domain,
\begin{equation}
\label{Eq:GenUnifBoundG}
	S
	:=
		\sup_{F\in\Fcal}\|\Gscr (F)\|_{L^\infty(\Dcal)}<\infty.
\end{equation}
As observed in \eqref{Eq:LipCondG}, the elliptic inverse problem considered in this paper falls in this general framework, which also encompasses other examples of nonlinear inverse problems such as those involving the Schr\"odinger equation considered in \cite{N17, NVW18}, for which the results in this section would apply as well. It also includes many linear inverse problems such as the classical Radon transform, see \cite{NVW18}.

%
%
%

%------------------------------------------------------------------------------------------------------------------------
%------------------------------------------------------------------------------------------------------------------------
\subsection{General contraction rates in Hellinger distance}
\label{SubSec:HellDistContr}
%------------------------------------------------------------------------------------------------------------------------
%------------------------------------------------------------------------------------------------------------------------

	Using the same notation as in Section \ref{Subsec:DivFormPDE}, and given a sequence of Borel prior probability measures $\Pi_N$ on $\Fcal$, we write $\Pi_N(\cdot|Y^{(N)},X^{(N)})$ for the posterior distribution of $F|(Y^{(N)},X^{(N)})$ (arising as after \eqref{Eq:PostDistr} and \eqref{Eq:JointLogLikelihood}).  We also continue to use the notation $p_F$ for the densities from \eqref{Eq:SingleLikelihood} now in the general observation model \eqref{Eq:GenDiscrObs} (and implicitly assume that the map $(F,(y,x)) \mapsto p_F(y,x)$ is jointly measurable to ensure existence of the posterior distribution). Below we formulate a general contraction theorem in the Hellinger distance  that forms the basis of the proofs of the main results. It closely follows the general theory in \cite{GvdV17} and its adaptation to the inverse problem setting in \cite{MNP19b} -- we include a proof for conciseness and convenience of the reader. 
	
	Define the Hellinger distance $h(\cdot, \cdot)$ on the set of probabilities density functions on $\R\times\Dcal$ (with respect to the product measure $dy\times dx$) by
$$
	h^2(p_1,p_2)
	:=
		\int_{\R\times\Dcal}\Big[\sqrt{p_1(y,x)}-\sqrt{p_2(y,x)}\Big]^2dydx.
$$
For any set $\Acal$ of such densities, let $N(\eta;\Acal,h), \ \eta>0,$ be the minimal number of Hellinger balls of radius $\eta$ needed to cover $\Acal$.

\begin{theorem}\label{Theo:GenHellRates}%------------------------------------------------------------------

	Let $\Pi_N$ be a sequence of prior Borel probability measures on $\Fcal$, and let $\Pi_N(\cdot |Y^{(N)},X^{(N)})$ be the resulting posterior distribution arising from observations $(Y^{(N)},X^{(N)})$ in model \eqref{Eq:GenDiscrObs}. Assume that for some fixed $F_0\in\Fcal$, and a sequence $ \delta_N >0$ such that $\delta_N\to0$ and $ \sqrt{N}\delta_N \to\infty$ as $N\to\infty$, the sets 
\begin{equation}
\label{Eq:Bn}
	\Bcal_N
	:=
		\Big\{F: \ E^1_{F_0}\Big[\log\frac{p_{F_0}(Y_1,X_1)}
		{p_F(Y_1,X_1)}\Big]\le \delta_N^2,
		\ E^1_{F_0}\Big[\log\frac{p_{F_0}(Y_1,X_1)}
		{p_F(Y_1,X_1)}\Big]^2\le \delta_N^2\Big\},
\end{equation}
satisfy for all $N$ large enough 
\begin{equation}
\label{Eq:Cond1SmallBall}
	\Pi_N(\Bcal_N )\ge a e^{-A N\delta_N^2},
	\quad \textrm{some $a,A>0$}.
\end{equation}
Further assume that there exists a sequence of Borel sets $\Acal_N\subset \Fcal$  for which
\begin{equation}
\label{Eq:Cond2ExcessMass}
	\Pi_N( \Acal_N^c)\lesssim e^{-B N\delta_N^2},
	\quad \textrm{some $B>A+2$}
\end{equation}
for all $N$ large enough, as well as
\begin{equation}
\label{Eq:Cond3MetricEntropy}
	\log N( \delta_N ; \Acal_N, h)\le C N\delta_N^2,
	\quad
	\textrm{some $C>0$.}
\end{equation}
Then, for sufficiently large $L=L(B,C)>4$ such that $L^2>12(B\vee C)$, and all $0<D<B-A-2$, as $N\to\infty$,
\begin{equation}
\label{Eq:GeneralPostContr}
	\Pi_N(F\in\Acal_N: h(p_F,p_{F_0})\le L \delta_N |Y^{(N)},X^{(N)})
	=
		1-O_{P^N_{F_0}}(e^{- D N\delta_N^2 }).
\end{equation}

\proof%-------------------------------------------------------------------------------------------------------------------

	We start noting that by Theorem 7.1.4 in \cite{GN16}, for each $L>4$ satisfying $L^2>12(B\vee C)$ we can find tests (random indicator functions) $\Psi_N=\Psi_N(Y^{(N)},X^{(N)})$ such that as $N\to\infty$
\begin{equation}
\label{Eq:ErrorBoundsMetricEntropyTests}
	E^N_{F_0}\Psi_N\to0, 
	\quad 	
	\sup_{F\in\Acal_N:h(p_F,p_{F_0})\ge L\delta_N}
	E^N_F(1-\Psi_N)\le e^{-B N\delta_N^2 }.
\end{equation}

	Next, denote the set whose posterior probability we want to lower bound by
$$
	\tilde\Acal_N=\{F\in\Acal_N :h(p_F,p_{F_0})\le L \delta_N \}; 
$$
and, using the first display in \eqref{Eq:ErrorBoundsMetricEntropyTests}, decompose the probability of interest as
\begin{align*}
	P^N_{F_0}&\Big(\Pi_N(\tilde\Acal_N^c|Y^{(N)},X^{(N)})\ge 
	e^{-D N\delta_N^2 }\Big)\\
	&=
		P^N_{F_0}\Big(\Pi_N(\tilde\Acal_N^c|Y^{(N)},X^{(N)})
		\ge e^{-DN\delta_N^2},\Psi_N=0\Big)\\
	&\quad
		+P^N_{F_0}\Big(\Pi_N(\tilde\Acal_N^c|Y^{(N)},X^{(N)})\ge 
		e^{-DN\delta_N^2},\Psi_N=1\Big)\\
	&=
		P^N_{F_0}\Big(\Pi_N(\tilde\Acal_N^c|Y^{(N)},X^{(N)})\ge e^{-DN\delta_N^2},\Psi_N=0\Big)
		+ o(1).
\end{align*}

	Next, let $\nu(\cdot)=\Pi_N(\cdot\cap \Bcal_N)/\Pi_N(\Bcal_N)$ be the restricted normalised prior on $\Bcal_N$, and define the event
\begin{equation}
\label{Eq:Cn}
	\Ccal_N=\Big\{\int_{\Bcal_N} \prod_{i=1}^N \frac{p_F}{p_{F_0}
	}(Y_i,X_i) d\nu(F)\ge e^{- 2 N\delta_N^2 }\Big\},
\end{equation}
for which Lemma 7.3.2 in \cite{GN16} implies that $P^N_{F_0}(\Ccal_N)\to1$ as $N\to\infty$. We then further decompose
\begin{align*}
	P^N_{F_0}\Big(\Pi_N(\tilde\Acal_N^c|Y^{(N)},X^{(N)})&
	\ge e^{-DN\delta_N^2},\Psi_N=0\Big)\\
	&=
		P^N_{F_0}\Big(\Pi_N(\tilde\Acal_N^c|Y^{(N)},X^{(N)})\ge e^{-DN\delta_N^2},
		\Psi_N=0,\Ccal_N\Big)+o(1)
\end{align*}
and in view of condition \eqref{Eq:Cond1SmallBall} and the above definition  of $\Ccal_N$, we see that
\begin{align*}
	P^N_{F_0}&\Big(\Pi_N(\tilde\Acal_N^c|Y^{(N)},X^{(N)})
	\ge e^{-DN\delta_N^2 },\Psi_N=0,\Ccal_N\Big)\\
	&=
		P^N_{F_0}\Bigg( \frac{\int_{\tilde\Acal_N^c}
		\prod_{i=1}^N p_F/p_{F_0}(Y_i,X_i)d\Pi_N(F)}
		{\int_{\Fcal}\prod_{i=1}^N p_F/p_{F_0}(Y_i,X_i)d\Pi_N(F)}
		\ge e^{-DN\delta_N ^2},\Psi_N=0,\Ccal_N\Bigg)\\
%	&=
%		P^N_{F_0}( \frac{\int_{\tilde\Acal_N^c}(1-\Psi_N)
%		\prod_{i=1}^N \frac{p_F}{p_{F_0}}(Y_i,X_i)d\Pi_N(F)}
%		{\Pi_N(\Bcal_N)\int_{\Fcal}\prod_{i=1}^N \frac{p_F}
%		{p_{F_0}}(Y_i,X_i)		\nu(dF)}
%		\ge e^{-DN\delta_N ^2},\Psi_N=0,\Ccal_N)\\
	&\le
		P^N_{F_0}\Bigg( \frac{\int_{\tilde\Acal_N^c}(1-\Psi_N)
		\prod_{i=1}^N p_F/p_{F_0}(Y_i,X_i)d\Pi_N(F)}
		{\int_{\Bcal_N}\prod_{i=1}^N p_F/p_{F_0}(Y_i,X_i)d\nu(F)}
		\ge \Pi_N(\Bcal_N)e^{-DN\delta_N ^2},\Ccal_N\Bigg)\\
%	&\le
%		P^N_{F_0}\Bigg(\int_{\tilde\Acal_N^c}(1-\Psi_N)
%		\prod_{i=1}^N \frac{p_F}{p_{F_0}}(Y_i,X_i)d\Pi_N(F)
%		\ge ae^{-(A+D+2)N\delta_N ^2},\Ccal_N\Bigg)\\
	&\le
		P^N_{F_0}\Bigg(\int_{\tilde\Acal_N^c}(1-\Psi_N)
		\prod_{i=1}^N \frac{p_F}{p_{F_0}}(Y_i,X_i)d\Pi_N(F)
		\ge ae^{-(A+D+2)N\delta_N ^2}\Bigg).
\end{align*}
 We conclude applying Markov's inequality and Fubini's theorem, jointly with the fact that for all $F\in\Fcal$
$$
	E^N_{F_0}\Big[(1-\Psi_N)\prod_{i=1}^N \frac{p_F}{p_{F_0}}
	(Y_i,X_i)\Big]
	=
		E^N_{F_0}\Big[(1-\Psi_N)\prod_{i=1}^N
		\frac{dP^{1}_F}{dP^{1}_{F_0}}(Y_i,X_i)\Big]
	=
		E^N_F[1-\Psi_N],
$$
to upper bound the last probability by
\begin{align*}
		\frac{1}{a}&e^{(A+D+2) N\delta_N^2 }\Big(\int_{\Acal_N^c}E^N_F
		[1-\Psi_N]d\Pi_N(F)
		+ \int_{\{F\in\Acal_N:h(p_{F_0},p_F)>L \delta_N \}}E^N_F
		[1-\Psi_N]d\Pi_N(F)\\
	&\ 
	 	+ \int_{\{F\in\Acal_N^c:h(p_{F_0},p_F)>L \delta_N \}}
		E^N_F[1-\Psi_N]d\Pi_N(F)\Big)\\
	&\  \le
		\frac{1}{a}e^{(A+D+2) N\delta_N^2 }\Big(2\Pi_N(\Acal_N^c)
	 	+ \int_{\{F\in\Acal_N:h(p_{F_0},p_F)>L \delta_N \}}
		E^N_F[1-\Psi_N]d\Pi_N(F)\Big)\\
	&\  \lesssim
		e^{-(B-A-D-2)N\delta_N^2}=o(1)
\end{align*}
as $N\to\infty$ since $B>A+D+2$, having used the excess mass condition \eqref{Eq:Cond2ExcessMass} and the second display in \eqref{Eq:ErrorBoundsMetricEntropyTests}.

\endproof%--------------------------------------------------------------------------------------------------------------

\end{theorem}%--------------------------------------------------------------------------------------------------------

%
%
%

%------------------------------------------------------------------------------------------------------------------------
%------------------------------------------------------------------------------------------------------------------------
\subsection{Contraction rates for rescaled Gaussian priors}
\label{SubSec:HellDistContr}
%------------------------------------------------------------------------------------------------------------------------
%------------------------------------------------------------------------------------------------------------------------

	While the previous result assumed a general sequence of priors, we now derive explicit contraction rates in $L^2-$prediction risk for the specific choices of priors considered in Section \ref{Sec:Results}. We start with the ‘re-scaled' priors of Section \ref{Subsec:GaussianPriors}.
	
\begin{theorem}\label{Theo:GenFwdRates1}%---------------------------------------------------------------

	Let the forward map $\Gscr$  satisfy \eqref{Eq:GenLipCondG} and \eqref{Eq:GenUnifBoundG} for given $\beta,\gamma,\kappa,\ge0$ and $S>0$. For integer $\alpha>\beta+d/2,$ consider a Gaussian prior $\Pi_N$ constructed as in \eqref{Eq:Prior1} with scaling $N^{d/(4\alpha+4\kappa+2d)}$ and with base prior $F'\sim\Pi'$ satisfying Condition \ref{Cond:BasePrior1}  with RKHS $\mathcal H$. Let $\Pi_N(\cdot|Y^{(N)},X^{(N)})$ be the resulting posterior arising from observations $(Y^{(N)},X^{(N)})$ in \eqref{Eq:GenDiscrObs}, assume $F_0\in\Hcal$ and set $\delta_N=N^{-(\alpha+\kappa)/(2\alpha+2\kappa+d)}$.
	
	Then for any $D>0$ there exists $L>0$ large enough (depending on $\sigma,F_0,D, \alpha,$ and $\beta,\gamma,\kappa,S,d$)  such that, as $N\to\infty$,
\begin{equation}
\label{Eq:GenFwdRates1}
	\Pi_N(F:\|\Gscr (F)-\Gscr (F_0)\|_{L^2(\Dcal)}>L\delta_N|Y^{(N)},X^{(N)}) 
	=
		O_{P^{N}_{F_0}}(e^{-	DN\delta_N^2}),
\end{equation}
and for sufficiently large $M>0$ (depending on $\sigma,D,\alpha,\beta,\gamma,\kappa,d$)
\begin{equation}
\label{Eq:GenRegularisation1}
	\Pi_N(F: \|F\|_{C^\beta}>M|Y^{(N)},X^{(N)})=O_{P^{N}_{F_0}}(e^{-DN\delta_N^2}).
\end{equation}

\end{theorem}%-------------------------------------------------------------------------------------------------------

\begin{remark}	\normalfont Inspection of the proof (cf. after \eqref{Eq:GenMetrEntr}) shows that if $\kappa=0$ in \eqref{Eq:GenLipCondG}, then the RKHS $\Hcal$ in Condition \ref{Cond:BasePrior1} can be assumed to be continuously imbedded in $H^\alpha(\Ocal)$ instead of $H^\alpha_c(\Ocal)$. The same remark in fact applies for $\kappa<1/2$.
\end{remark}

\proof%------------------------------------------------------------------------------------------------------------------

	In view of the boundedness assumption \eqref{Eq:GenUnifBoundG} on $\Gscr$, we have by Lemma \ref{Lemma:HellingerAndL2SmallBalls} below that for some $q>0$ (depending on $\sigma,S$)
\begin{align*}
	E^{1}_{F_0}\Big[
	\log\frac{p_{F_0}(Y_1,X_1)}{p_F(Y_1,X_1)}\Big]\vee
	E^{1}_{F_0}\Big[\log\frac{p_{F_0}(Y_1,X_1)}
	{p_F(Y_1,X_1)}\Big]^2
	\le 
		q\|\Gscr (F_0)-\Gscr (F)\|^2_{L^2(\Dcal)}.
\end{align*}
Hence, for $\Bcal_N$ the sets from \eqref{Eq:Bn} we have
$
	\{F : \|\Gscr (F_0)-\Gscr (F)\|_{L^2(\Dcal)}
	\le \delta_N/q\}
	\subseteq
		 \Bcal_N,
$
which in turn implies the small ball condition \eqref{Eq:Cond1SmallBall} since by Lemma  \ref{Lemma:SmallBall1} below (premultiplying, if needed, $\delta_N$ by a sufficiently large but fixed constant)
$$
	\Pi_N(F: \|\Gscr (F)-\Gscr (F_0)\|_{L^2(\Dcal)} \le\delta_N/q)
	\gtrsim
		e^{-AN\delta_N^2}
$$
for some $A>0$ and all $N$ large enough. Next, for all $D>0$ and any $B>A+D+2$, we can choose sets $\Acal_N$ as in Lemmas \ref{Lemma:ApproxSets1} and \ref{Lemma:MetricEntropy1} and verify the excess mass condition \eqref{Eq:Cond2ExcessMass} as well as the complexity bound \eqref{Eq:Cond3MetricEntropy}. Note that $\|F\|_{C^\beta}\le M$ for all $F\in\Acal_N$. We then conclude by Theorem \ref{Theo:GenHellRates} that for some $L'>0$ large enough
$$
\Pi_N(F\in\Acal_N: h(p_F,p_{F_0})\le L' \delta_N |Y^{(N)},X^{(N)})
	=
		1-O_{P^N_{F_0}}(e^{- D N\delta_N^2 })
$$
yielding the claim for some appropriate $L>0$ using the first inequality in \eqref{Eq:HellL2}.

\endproof%-------------------------------------------------------------------------------------------------------------

	The following key lemma shows that the (non-Gaussian) prior induced on the regression functions $\mathscr G(F)$ assigns sufficient mass to a $L^2$-neighbourhood of $\mathscr G(F_0)$.

\begin{lemma}\label{Lemma:SmallBall1}%---------------------------------------------------------------------

	Let $\Pi_N,F_0$ and $\delta_N$ be as in Theorem \ref{Theo:GenFwdRates1}. Then, for sufficiently large $c>0$  there exists $A>0$ (depending on $c,F_0, \alpha, \beta,\gamma,\kappa,d$) such that
\begin{equation}
\label{Eq:SmallBall1}
	\Pi_N(F: \|\Gscr (F)-\Gscr (F_0)\|_{L^2(\Dcal)} \le c \delta_N)
	\gtrsim
		e^{-AN\delta_N^2}
\end{equation}
for all $N$ large enough.

\proof%------------------------------------------------------------------------------------------------------------------

	Since $F_0\in\Hcal$, $\|F_0\|_{C^\beta}<\infty$ by a Sobolev imbedding. Let  $M>\|F_0\|_{C^\beta}\vee 1$ be a fixed constant. Using \eqref{Eq:GenLipCondG}, we obtain for some $k>0$
\begin{equation*}
\begin{split}
	\Pi_N(F: \|\Gscr (F)-\Gscr &(F_0)\|_{L^2(\Dcal)} \le  c \delta_N)\\
	&\ge 
	\Pi_N(F: \|F-F_0\|_{(H^\kappa)^*} \le  c k M^{-\gamma} \delta_N, 
	\|F-F_0\|_{C^\beta} \le M )\\
	&=
		\Pi_N(F:F-F_0\in C_1\cap C_2),
\end{split}
\end{equation*}
where
$$
	C_1:=\{F: \|F\|_{(H^\kappa)^*}\le  c kM^{-\gamma} \delta_N\},
	\quad 
	C_2:=\{F: \|F\|_{C^\beta}\le M\}.
$$
Then, recalling that the RKHS $\Hcal_N$ of $\Pi_N$ coincides with $\Hcal$ with RKHS norm $\|\cdot\|_{\Hcal_N}$ given in \eqref{Eq:RKHS1}, now with scaling $N^{d/(4\alpha+4\kappa+2d)} = \sqrt N \delta_N$, we can use Corollary 2.6.18 in \cite{GN16} to lower bound the last probability by
\begin{align*}
	e^{-\|F_0\|^2_{\Hcal_N}/2}	\Pi_N( C_1\cap C_2)\nonumber
	&=
		e^{-\frac{1}{2}N \delta_N^2\|F_0\|^2_{\Hcal}} \Pi_N( C_1\cap C_2) \\
		& \ge e^{-\frac{1}{2}N \delta_N^2\|F_0\|^2_{\Hcal}}
		 \big(\Pi_N(C_1) - \Pi_N(C^c_2) \big)
\end{align*}

	To upper bound $\Pi_N(C^c_2)$, note that by construction of $\Pi_N$\begin{align*}
	\Pi_N( C_2^c)
	&=
		\Pr( \| F'\|_{C^\beta}>  M N\delta_N^2),
	\quad F'\sim \Pi'.
\end{align*}
By Condition \ref{Cond:BasePrior1}, $F'$ defines a centred Gaussian Borel random element in a separable measurable subspace $\Ccal$ of $C^\beta$, and by the Hahn-Banach theorem and the separability of $\Ccal$, $\|F'\|_{C^\beta}$ can then be represented as a countable supremum
$$
	\| F'\|_{C^\beta}=\sup_{T\in\Tcal}|T(F')|
$$
of actions of bounded linear functionals $\Tcal=(T_m)_{m\in\N}\subset (C^\beta)^*$.
It follows that the collection $\{T_m(F')\}_{m\in\N}$ is a centred Gaussian process with almost surely finite supremum, so that by Fernique's theorem \cite[Theorem 2.1.20]{GN16}
$$
	E\|F'\|_{C^\beta}=E\sup_{m\in\N}|T_m( F')|<\infty;
	\quad
	\tau^2:=\sup_{m\in\N}E|T_m(F')|^2<\infty.
$$
We then apply the Borell-Sudakov-Tirelson inequality \cite[Theorem 2.5.8]{GN16} to obtain for all $N$ large enough, 
\begin{align}
\label{Eq:BTSIneq}
	\Pr\big(\|F'\|_{C^\beta}\ge  M\sqrt{N}\delta_N\big)
	&\le
		\Pr\big(\|F'\|_{C^\beta}\ge E\|F'\|_{C^\beta} 
		+ M\sqrt{N}\delta_N/2\big)
	\le
		e^{-\frac{1}{8}(M/\tau)^2N\delta_N^2}.
\end{align}

	We proceed finding a lower bound for the prior probability of $C_1$, which, by construction of $\Pi_N$, satisfies
\begin{align*}
		\Pi_N(F\in C_1)
	=
		\Pi'( F':\| F'\|_{(H^\kappa)^*}\le c k M^{-\gamma} \sqrt{N}\delta_N^2).
\end{align*}
For any integer $\alpha>0$ and any $\kappa\ge0$, letting $B^\alpha_c(r):=\{F\in H^\alpha_c, \ \|F\|_{H^\alpha}\le r\}, \ r>0,$ we have the metric entropy estimate:
\begin{equation}
\label{Eq:GenMetrEntr}
	\log N(\eta; B_c^\alpha(r),\|\cdot\|_{(H^\kappa)^*})
	\lesssim
		 (r/\eta)^{d/(\alpha+\kappa)} \quad \forall \eta>0;
\end{equation}
see the proof of Lemma 19 in \cite{NVW18} for the case $\kappa\ge1/2$, and Theorem
4.10.3 in \cite{Triebel1978} for $\kappa<1/2$ (in the latter case, we note in fact that the estimate holds true also for balls in the whole space $H^\alpha$). Hence, since $\Hcal$ is continuously imbedded into $H^\alpha_c$, letting $B_\Hcal(1)$ be the unit ball of $\Hcal$, we have $B_\Hcal(1)\subseteq B_c^\alpha(r)$ for some $r>0$, implying that for all $\eta>0$
\begin{equation}
\label{Eq:MetricEntropyEstim}
	\log N(\eta; B_\Hcal(1),\|\cdot\|_{(H^\kappa)^*})
	\le 
		\log N(\eta; B_c^\alpha(r),\|\cdot\|_{(H^\kappa)^*})
	\lesssim
 		\eta^{-d/(\alpha+\kappa)}.
\end{equation}
Then, for all $N$ large enough, the small ball estimate in Theorem 1.2 in \cite{LL99} yields 
\begin{align*}
	-\log\Pi'(& F':\| F'\|_{(H^\kappa)^*}\le c k M^{-\gamma} \sqrt{N}\delta_N^2)
	\nonumber\\
	&\lesssim
		( c M^{-\gamma} \sqrt{N}\delta_N^2)^{-2\frac{d}{\alpha+\kappa}
		(2-d/(\alpha+\kappa))^{-1}}\nonumber\\
	&= 
		(M^\gamma/c)^\frac{2d}{2\alpha+2\kappa-d}
		[N^{-(\alpha+\kappa-d/2)/(2\alpha+2\kappa+d)}]
		^{-\frac{2d}{2\alpha+2\kappa-d}}\nonumber\\
	&=
		(M^\gamma/c)^\frac{2d}{2\alpha+2\kappa-d}
		N\delta_N^2.
\end{align*}

Thus, for $k'>0$ a fixed constant, we obtain the lower bound
\begin{align*}
\Pi_N(F&: \|\Gscr (F)-\Gscr (F_0)\|_{L^2(\Dcal)} \le  c \delta_N)\\
&\ge
 	e^{-\frac{1}{2}N \delta_N^2\|F_0\|^2_{\Hcal}}
	\Big(e^{-k' (M^\gamma/c)
	^\frac{2d}{2\alpha+2\kappa-d}N\delta_N^2}
	- e^{-\frac{1}{8}(M/\tau)^2N\delta_N^2}\Big)\\
&\gtrsim
	e^{-A N\delta_N^2}
\end{align*}
having taken $c>0$ large enough (satisfying
$k' (M^\gamma/c)^\frac{2d}{2\alpha+2\kappa-d}<\frac{1}{8}(M/\tau)^2$), and where 
$A= \frac{1}{2}\|F_0\|^2_{\Hcal} + k' (M^\gamma/c)^\frac{2d}{2\alpha+2\kappa-d}$.

\endproof%--------------------------------------------------------------------------------------------------------------

\end{lemma}%----------------------------------------------------------------------------------------------------------

	We now construct suitable approximating sets for which we check the excess mass condition \eqref{Eq:Cond2ExcessMass}.
	
\begin{lemma}\label{Lemma:ApproxSets1}%---------------------------------------------------------------

	Let $\Pi_N$ and $\delta_N$ be as in Theorem \ref{Theo:GenFwdRates1}. Define for any $M,Q>0$ 
\begin{equation}
\label{Eq:ApproxSets1}
	\Acal_N
	= 
		\{F:F= F_1 +F_2: \|F_1\|_{(H^\kappa)^*} \le Q   \delta_N,\
		\|F_2\|_{\Hcal} \le M,\ \|F\|_{C^\beta} \le M\}.
\end{equation}
Then for any $B>0$ and for sufficiently large $M, Q$ (both depending on $B,\alpha,\beta,\gamma,\kappa,d$), for all $N$ large enough,
\begin{equation}
\label{Eq:ExpDecay}
	\Pi_N(\Acal^c_N) \le 2e^{-B N\delta_N^2}.
\end{equation}

\proof%------------------------------------------------------------------------------------------------------------------

	By \eqref{Eq:BTSIneq}, taking $M\gtrsim \sqrt{ B}$, we obtain for all $N$ large enough that
$
	\Pi_N(F: \|F\|_{C^\beta}\le M)\ge 1-e^{-BN\delta^2_N}.
$
Thus, the claim will follow if we can derive a similar lower bound for
\begin{align*}
	\Pi_N(F:F&= F_1 +F_2: \|F_1\|_{(H^\kappa)^*} \le Q \delta_N, \|F_2\|_{\Hcal} \le M  )\\
	&=
		 \Pi' (F':F' =  F'_1 + F'_2, \ \|F'_1\|_{(H^\kappa)^*} \le Q \sqrt{N}\delta_N^2,
		  \|F'_2\|_\Hcal \le M\sqrt{N}\delta_N),
\end{align*}
having used that $N^{d/(4\alpha+4\kappa+d)}=\sqrt{N}\delta_N$. Using Theorem 1.2 in \cite{LL99} as after \eqref{Eq:MetricEntropyEstim}, we deduce that for some $q>0$
$$
	-\log\Pi'(F':\|F'\|_{(H^\kappa)^*}\le Q \sqrt{N}\delta_N^2)
	\le q
	( Q \sqrt{N}\delta_N^2)^{-\frac{2d}{2\alpha+2\kappa-d}}
$$
so that for any $Q>(B/q)^{-(2\alpha+2\kappa-d)/(2d)}$
\begin{equation}
\label{Eq:ExpLB1}
	-\log\Pi'(F':\|F'\|_{(H^\kappa)^*}\le Q \sqrt{N}\delta_N^2)
	\le 
		B( \sqrt{N}\delta_N^2)^{-\frac{2d}{2\alpha+2\kappa-d}}
%	&=
%		[N^{-(\alpha+1-d/2)/(2\alpha+2+d})]^{-\frac{2d}{2\alpha+2\kappa-d}}\\
%	&=
%		BN^{\frac{d}{2\alpha+2\kappa+d}}\\
	=
		BN\delta_N^2.
\end{equation}
Next, denote
$$
	M_N=-2\Phi^{-1}(e^{- B N\delta_N^2})
$$
where $\Phi$ is the standard normal cumulative distribution function. Then by standard inequalities for $\Phi^{-1}$ we have $M_N\simeq \sqrt{BN}\delta_N$ as $N\to\infty$, so that taking $M\gtrsim\sqrt{B}$ implies
\begin{align*}
	 \Pi'(F':F' &=  F'_1 + F'_2, \ \|F'_1\|_{(H^\kappa)^*} \le Q \sqrt{N}\delta_N^2,
	  \|F'_2\|_\Hcal \le M\sqrt{N}\delta_N)\\
	&\ge
		\Pi' (F':F' =  F'_1 + F'_2,\ \|F'_1\|_{(H^\kappa)^*} \le Q\sqrt{N}\delta_N^2,
		 \|F'_2\|_\Hcal \le M_N).
\end{align*}
By the isoperimetric inequality for Gaussian processes \cite[Theorem 2.6.12]{GN16}, the last probability is then lower bounded, using  \eqref{Eq:ExpLB1}, by
\begin{align*}
	\Phi(\Phi^{-1}[\Pi'(F':\|F'\|_{(H^\kappa)^*}\le Q\sqrt{N}\delta_N^2)]+M_N)
	&\ge
		\Phi(\Phi^{-1}[e^{- B N\delta_N^2}]+M_N) = 1-e^{- B N\delta_N^2},
\end{align*}
concluding the proof.
		
\endproof%--------------------------------------------------------------------------------------------------------------

\end{lemma}%----------------------------------------------------------------------------------------------------------

	We conclude with the verification of the complexity bound \eqref{Eq:Cond3MetricEntropy} for the sets $\Acal_N$.
		
\begin{lemma}\label{Lemma:MetricEntropy1}%---------------------------------------------------------------

	Let $\Acal_N$ be as in Lemma \ref{Lemma:ApproxSets1} for some fixed $M,Q>0$. Then, 
$$
	\log N( \delta_N;\Acal_N,h)\le C N\delta_N^2,
$$
for some constant $C>0$ (depending on $\sigma, M,Q,\alpha,\beta,\gamma,\kappa,d,S$) and all $N$ large enough.

\proof%-------------------------------------------------------------------------------------------------------------------

	If $F\in \Acal_N$, then $F=F_1+F_2$ with 
$
	\|F_1\|_{(H^\kappa)^*}
	\le
		 Q \delta_N
$
and
$
		\|F_2\|_{H^\alpha}
		\le 
			M',
$
the latter inequality following from the continuous imbedding of $\Hcal$ into $H^\alpha_c$. Thus, recalling the metric entropy estimate \eqref{Eq:GenMetrEntr}, if 
 $$
 	\{H_1,\dots,H_P\}\subset B^\alpha_c(M'), 
	\quad 
	P\le e^{-q\delta_N^{-d/(\alpha+\kappa)}}=e^{-qN\delta_N^2},~q>0,
$$ 
 is a $\delta_N$-net with respect to $\|\cdot\|_{(H^\kappa)^*}$, we can find $H_i$ such that $\|F_2-H_i\|_{(H^\kappa)^*}\le\delta_N$. Then, using the second inequality in \eqref{Eq:HellL2} below and the local Lipschitz estimate \eqref{Eq:GenLipCondG},
\begin{align*}
	h(p_F,H_i)
	&\lesssim 
		\|\Gscr(F)-\Gscr(H_i)\|_{L^2(\Dcal)}\\
	&\lesssim
		(1+\|F\|^\gamma_{C^\beta}\vee\|H_i\|^\gamma_{C^\beta})\|F-H_i\|_{(H^\kappa)^*}.
\end{align*}
Recalling that if $F\in\Acal_N$ then also $\|F\|_{C^\beta}\le M$,  and using the Sobolev imbedding of $H^\alpha$ into $C^\beta$ to bound $\|H_i\|_{C^\beta}$, we then obtain				
\begin{align*}
	h(p_F,H_i)
	\lesssim 
	\|F-H_i\|_{(H^\kappa)^*}
	\lesssim 
		\|F-F_2\|_{(H^\kappa)^*}+\|F_2-H_i\|_{(H^\kappa)^*}
	\lesssim
		\delta_N.
\end{align*}
It follows that $\{H_1,\dots,H_P\}$ also forms a $q'\delta_N$-net for $\Acal_N$ in the Hellinger distance for some $q'>0$, so that
$$
	\log N( \delta_N;\Acal_N,h)
	\le
		\log N( \delta_N/q';B^\alpha_c(M),\|\cdot\|_{(H^\kappa)^*})
	\lesssim 
		 N\delta_N^2.
$$

\endproof%--------------------------------------------------------------------------------------------------------------

\end{lemma}%----------------------------------------------------------------------------------------------------------

%
%
%

%------------------------------------------------------------------------------------------------------------------------
%------------------------------------------------------------------------------------------------------------------------
\subsection{Contraction rates for hierarchical Gaussian series priors}
\label{SubSec:HellDistContr}
%------------------------------------------------------------------------------------------------------------------------
%------------------------------------------------------------------------------------------------------------------------

	We now derive contraction rates in $L^2$-prediction risk in the inverse problem \eqref{Eq:GenDiscrObs}, for the truncated Gaussian random series priors introduced in Section \ref{Subsec:RandSeriesPriors}. The proof again proceeds by an application of Theorem \ref{Theo:GenHellRates}.

\begin{theorem}\label{Theo:GenFwdRates2}%---------------------------------------------------------------

Let the forward map $\Gscr$  satisfy \eqref{Eq:GenLipCondG} and \eqref{Eq:GenUnifBoundG} for given $\beta,\gamma,\kappa\ge0$ and $S>0$. For any $\alpha>\beta+d/2$, let $\Pi$ be the random series prior in \eqref{Eq:Prior2}, and let $\Pi(\cdot|Y^{(N)},X^{(N)})$ be the resulting posterior distribution arising from observations $(Y^{(N)},X^{(N)})$ in \eqref{Eq:GenDiscrObs}. Then, for each $\alpha_0\ge\alpha$, any $F_0\in H^{\alpha_0}_K(\Ocal)$ and any $D>0$ there exists $L>0$ large enough (depending on $\sigma,F_0,D,\alpha,\beta,\gamma,\kappa,S,d$) such that, as $N \to \infty$,
\begin{equation}
\label{Eq:GenFwdRates1}
	\Pi(F:\|\Gscr (F)-\Gscr (F_0)\|_{L^2(\Dcal)}>L\xi_N|Y^{(N)},X^{(N)})
	=
		O_{P^{N}_{F_0}}(e^{-DN\xi_N^2}),
\end{equation}
where $\xi_N= N^{-(\alpha_0+\kappa)/(2\alpha_0+2\kappa+d)}\log N$. Moreover, for $\Hcal_J$ the finite-dimensional subspaces from \eqref{Eq:SubspaceHj} and $J_N\in\N$ such that $2^{J_N}\simeq N^{1/(2\alpha_0+2\kappa+d)}$, we also have that for sufficiently large $M>0$ (depending on $D,\alpha,\beta,d$)
\begin{equation}
\label{Eq:Regularisation2}
	\Pi(F: F\in\Hcal_{J_N}, \ \|F\|_{H^\alpha}\le M 2^{J_N\alpha}N\xi_N^2 |Y^{(N)},X^{(N)})=O_{P^{N}_{F_0}}(e^{-DN\xi_N^2}).
\end{equation}

\end{theorem}%--------------------------------------------------------------------------------------------------------

	We begin deriving a suitable small ball estimate in the $L^2$-prediction risk.
	
\begin{lemma}\label{Lemma:SmallBall2}%------------------------------------------------------------------
	Let $\Pi,\ F_0$ and $\xi_N$ be as in Theorem \ref{Theo:GenFwdRates2}. Then, for sufficiently large $q>0$ there exists $A>0$ (depending on $q,F_0,\alpha,\beta,\gamma,\kappa,d)$ such that
\begin{equation}
\label{Eq:SmallBall2}
	\Pi(F: \|\Gscr (F)-\Gscr (F_0)\|_{L^2(\Dcal)} \le q\xi_N)
	\gtrsim
		e^{-AN\xi_N^2}
\end{equation}
for all $N$ large enough.

\proof%-----------------------------------------------------------------------------------------------------------------

	For each  $j\in\N$, denote by $\Pi_j$ the Gaussian probability measure on the finite dimensional subspace $\Hcal_j$ in \eqref{Eq:SubspaceHj} defined as after \eqref{Eq:Prior2} with the series truncated at $j$. For  $J_N\in\N:2^{J_N}\simeq N^{1/(2\alpha_0+2\kappa+d)}$, note
\begin{equation}
\label{Eq:NXiN}
	2^{J_Nd}\log 2^{J_Nd}
	\simeq
		 N^{d/(2\alpha_0+2\kappa+d)}\log N=N\xi_N^2,
\end{equation}
so that, recalling the properties \eqref{Eq:PropertiesOfJ} of the random truncation level $J$, for some $s>0,$
$$
	\Pr(K=J_N)\gtrsim e^{-2^{J_Nd}\log 2^{J_Nd}} \ge e^{-sN\xi_N^2}
$$
for all $N$ large enough. It follows
\begin{align*}
	\Pi(F:\|\Gscr (F)-\Gscr (F_0)\|_{L^2}\le  q\xi_N)
%	&=
%		\sum_{j=1}^\infty\Pi_k(F:\|\Gscr (F)-\Gscr (F_0)\|_{L^2}\le \xi_N)\Pr(J=j)\\
	&\ge
		\Pi_{J_N}(F:\|\Gscr (F)-\Gscr (F_0)\|_{L^2}\le q\xi_N)\Pr(K=J_N)\\
%	&\gtrsim
%		\Pi_{J_N}(F:\|\Gscr (F)-\Gscr (F_0)\|_{L^2}\le \xi_N/q)e^{-2^{J_Nd}\log 2^{J_Nd}}\\
	&\gtrsim
		\Pi_{J_N}(F:\|\Gscr (F)-\Gscr (F_0)\|_{L^2}\le q\xi_N)e^{-sN\xi_N^2}.
\end{align*}

	Next, let 
$$
	P_{\Hcal_{J_N}}(F_0)=\chi\sum_{\ell\le {J_N},r\in\Rcal_\ell}
	\langle F_0,\Psi_{\ell r}\rangle_{L^2}\Psi_{\ell r}
$$
 be the `projection' of $F_0$ onto $\Hcal_{J_N}$. Since $F_0\in H^{\alpha_0}_K \subset C^\beta$ by a Sobolev imbedding, it follows using \eqref{Eq:GenLipCondG} and standard approximation properties of wavelets (cf. \eqref{Eq:DualNormApprox}), 
$$
	\|\Gscr (F_0)-\Gscr (P_{\Hcal_{J_N}}(F_0))\|_{L^2(\Dcal)}
	\lesssim
		\|F_0-P_{\Hcal_{J_N}}(F_0)\|_{(H^\kappa)^*}
	\lesssim 
		2^{-J_N(\alpha_0+\kappa)}
	=
		N^{-\frac{\alpha_0+\kappa}{2\alpha_0+2\kappa+d}},
$$
which implies by the triangle inequality that
\begin{align*}
	\Pi_{J_N}&(F:\|\Gscr (F)-\Gscr (F_0)\|_{L^2}\le q \xi_N)\\
	&\ge
		\Pi_{J_N}(F:\|\Gscr (F)-\Gscr(P_{\Hcal_{J_N}}(F_0))\|_{L^2}
		\le q\xi_N-\|\Gscr (F_0)-\Gscr(P_{\Hcal_{J_N}}(F_0))\|_{L^2})\\
	&\ge
		\Pi_{J_N}(F:\|\Gscr (F)-\Gscr(P_{\Hcal_{J_N}}(F_0))\|_{L^2}\le q'\xi_N).
\end{align*}
Using again that $H^{\alpha}$ imbeds continuously into $C^\beta$ as well as \eqref{Eq:GenLipCondG} and \eqref{Eq:AnotherBound}, we can lower bound the last probability by
\begin{align*}
	\Pi_{J_N}(F&:\|\Gscr (F)-\Gscr(P_{\Hcal_{J_N}}(F_0))\|_{L^2(\Dcal)}
	\le q'\xi_N,\|F-P_{\Hcal_{J_N}}(F_0)\|_{H^\alpha(\Ocal)}\le \xi_N)\\
	&\ge
		\Pi_{J_N}(F:\|F-P_{\Hcal_{J_N}}(F_0)\|_{(H^\kappa(\Ocal))^*}\le q''\xi_N,
		\|F-P_{\Hcal_{J_N}}(F_0)\|_{H^\alpha(\Ocal)}\le \xi_N)\\
	&\ge
		\Pi_{J_N}(F:\|F-P_{\Hcal_{J_N}}(F_0)\|_{H^\alpha(\Ocal)}\le q'''\xi_N),
\end{align*}
which, by Corollary 2.6.18 in \cite{GN16} and in view of \eqref{Eq:FiniteRKHSNorm} is further lower bounded by
$$
		e^{-\frac{1}{2}\|P_{\Hcal_{J_N}}(F_0)\|^2_{\Hcal_{J_N}}}	
		\Pi_{J_N}	(F:\|F\|_{H^\alpha}\le q'''\xi_N)\\
	\ge 
		e^{-s'\|F_0\|^2_{H^{\alpha_0}}}
		\Pi_{J_N}(F:\|F\|_{H^\alpha}\le q'''\xi_N).
$$
Now since $f \mapsto \chi f, \chi \in C^\infty(\Ocal)$, is continuous on $H^\alpha(\Ocal)$,
\begin{align*}
	\Pi_{J_N}(F:\|F\|_{H^\alpha}\le q'''\xi_N)
	&=
		\Pr\Big(\Big\|\chi \sum_{\ell\le J_N, r\in\Rcal_\ell} 
		2^{-\ell\alpha}F_{\ell r}\Psi_{\ell r}\Big\|_{H^\alpha}\le q'''\xi_N\Big)\\
	&\ge
		\Pr\Bigg(\sum_{m=1}^{\textnormal{dim}(\Hcal_{J_N})} Z_m^2\le t\xi_N^2\Bigg)
\end{align*}
for some $t>0$, where $Z_m\iid N(0,1)$, and where we have used the wavelet characterisation of the $H^\alpha(\mathbb R^d)$ norm. To conclude, note that the last probability is greater than
\begin{align*}
		\Pr&\Big(\sqrt{\textnormal{dim}(\Hcal_{J_N})}\max_{m\le \textnormal{dim}(\Hcal_{J_N})}
		 |Z_m| \le \sqrt t\xi_N\Big)\\
	&\ge
		\Pr\Big(\max_{m\le \textnormal{dim}(\Hcal_{J_N})} |Z_m| 
		\le t'N^{-\frac{\alpha_0+\kappa}{2\alpha_0+2\kappa+d}}
		N^{-\frac{d/2}{2\alpha_0+2\kappa+d}}\Big)	\\
	&=
		\prod_{m\le \textnormal{dim}(\Hcal_{J_N})}
		\Pr\Big(|Z_m| \le t'N^{-\frac{\alpha_0+\kappa+d/2}{2\alpha_0+2\kappa+d}}\Big).
\end{align*}
 Finally, a standard calculation shows that $\Pr(|Z_1| \le t)\gtrsim t$ if $t\to0$, and hence the last product is lower bounded, for large $N$, by
 \begin{align*}
		\left(t'N^{-\frac{\alpha_0+\kappa+d/2}{2\alpha_0+2\kappa+d}}
		\right)^{\textnormal{dim}(\Hcal_{J_N})}
	=
		e^{\textnormal{dim}(\Hcal_{J_N})\log\big( t'N^{-\frac{\alpha_0+\kappa+d/2}
		{2\alpha_0+2+d}}\big)}				
	\ge
		e^{-t''2^{J_Nd}\log N}
	=
		e^{-t'''N\xi_N^2}.
\end{align*}

\endproof%------------------------------------------------------------------------------------------------------------

\end{lemma}%--------------------------------------------------------------------------------------------------------

	In the following lemma we construct suitable approximating sets, for which we check the excess mass condition \eqref{Eq:Cond2ExcessMass} and the complexity bound \eqref{Eq:Cond3MetricEntropy} required in Theorem \ref{Theo:GenHellRates}.

\begin{lemma}\label{Lemma:ApproxSets2}%-------------------------------------------------------

	Let $\Pi, \ \xi_N$ and $J_N$ be as in Theorem \ref{Theo:FwdRates2}, and let $\Hcal_{J_N}$ be the finite dimensional subspace defined in \eqref{Eq:SubspaceHj} with $J=J_N$. Define for each $M>0$
\begin{equation}
\label{Eq:ApproxSets2}
	\Acal_N=\Big\{F\in \Hcal_{J_N}, \ \|F\|_{H^\alpha}\le M 2^{J_N\alpha}N\xi_N^2
	\Big\}.
\end{equation}
Then, for any $B>0$  there exists $M>0$ large enough (depending on $B,\alpha,\beta,d$) such that, for sufficiently large $N$
\begin{equation}
\label{Eq:ExcessMass2}
	\Pi(\Acal_N^c)\le 2e^{-BN\xi_N^2}. 
\end{equation}
Moreover, for each fixed $M>0$ and all $N$ large enough 
\begin{equation}
\label{Eq:MetricEntropy2}
\log N(\xi_N;\Acal_N,h)\le C N\xi_N^2
\end{equation}
for some $C>0$ (depending on $\sigma,\alpha,\beta,\gamma,\kappa,S,d$).

\proof%-------------------------------------------------------------------------------------------------------------------

	Letting $Z_m\iid N(0,1)$, noting $\|F\|^2_{H^\alpha}\le 2^{2J_N\alpha}\sum_{\ell \le J_N, r\in\Rcal_\ell} F_{\ell r}^2$ for all $F\in\Hcal_{J_N}$ (cf. \eqref{Eq:FiniteRKHSNorm}) and using \eqref{Eq:NXiN} and \eqref{Eq:PropertiesOfJ}, we have for sufficiently large $N$
\begin{align*}
	\Pi(\Acal_N^c)
	&\le
		\Pr(J> J_N)+
		\Pr\Big(\sum_{\ell\le J\land J_N, r\in\Rcal_\ell} F^2_{\ell r}\le MN\xi_N^2\Big)\\
	&\le
		e^{-2^{J_Nd}\log 2^{J_Nd}}
		+\Pr\Big(\sum_{m\le \textnormal{dim}(\Hcal_{J_N})}  Z_m^2 >MN\xi^2_N\Big)\\
	&\le 
		e^{-BN\xi_N^2}
		+\Pr\Big(\sum_{m\le \textnormal{dim}(\Hcal_{J_N})}  (Z_m^2-1) >\bar MN\xi^2_N\Big)
\end{align*}
for any constant $0<\bar M<M^2-1$, since $ \textnormal{dim}(\Hcal_{J_N})\lesssim 2^{J_Nd}\simeq N^{d/(2\alpha+2+d)}=o(N\xi_N^2)$. The bound \eqref{Eq:ExcessMass2} then follows applying Theorem 3.1.9 in \cite{GN16} to upper bound the last probability, for any $B>$ and for sufficiently large  $M$ and $\bar M$, by
$$
		e^{-\frac{\bar M^2(N\xi_N^2)^2}{4\textnormal{dim}(\Hcal_{J_N})+\bar MN\xi_N^2}}
%	\le 
%		e^{-\frac{(M')^2(N\xi_N^2)^2}{M''N\xi_N^2}}
	\le 
		e^{-BN\xi_N^2}.
$$

	We proceed with the derivation of \eqref{Eq:MetricEntropy2}. By choice of $J_N$, if $F\in\Acal_N$ then
$
	\|F\|_{H^\alpha}^2
	\lesssim
		N^{(2\alpha)/(2\alpha+2\kappa+d)} N\xi_N^2.
$
 Hence, by the second inequality in \eqref{Eq:HellL2}, using \eqref{Eq:GenLipCondG} and the Sobolev imbedding of $H^\alpha$ into $C^\beta$, if $F_1,F_2\in\Acal_N$ then
\begin{align*}
	h(p_{F_1},p_{F_2})
	&\lesssim
		\|\Gscr(F_1)-\Gscr(F_2)\|_{L^2(\Dcal)}\\
	&\lesssim
		(1+(N^\frac{\alpha}{2\alpha+2\kappa+d}\sqrt{N}\xi_N)^\gamma)
		\|F_1-F_2\|_{(H^\kappa)^*}\\
		&\lesssim
		N^\frac{\alpha \gamma}{2\alpha+2\kappa+d}(\sqrt{N}\xi_N)^\gamma
		\sqrt{\sum_{\ell\le J_N,r\in\Rcal_\ell} (F_{1,\ell r}- F_{2, \ell r})^2}.
\end{align*}
Therefore, using the standard metric entropy estimate for balls $B_{\R^p}(r), \ r>0,$ in Euclidean spaces  \cite[Proposition 4.3.34]{GN16}, we see that for $N$ large enough
\begin{align*}
	\log N(\xi_N;\Acal_N,h)
	&\lesssim	\log N\Big(\xi_NN^\frac{-\alpha \gamma}{2\alpha+2\kappa+d}
		(\sqrt{N}\xi_N)^{-\gamma}; 
		B_{\R^{\textnormal{dim}(\Hcal_{J_N})}}(M\sqrt{N}\xi_N),
		\|\cdot\|_{\R^{\textnormal{dim}(\Hcal_{J_N})}}\Big)\\
	&\le
		\textnormal{dim}(\Hcal_{J_N})
		\log \frac{3 M\sqrt{N}\xi_N}
		{\xi_NN^\frac{\alpha \gamma}{2\alpha+2\kappa+d}(\sqrt{N}\xi_N)^{-\gamma}} \lesssim N\xi_N^2.
\end{align*}

\endproof%-------------------------------------------------------------------------------------------------------------

\end{lemma}%--------------------------------------------------------------------------------------------------------

%
%
%

%------------------------------------------------------------------------------------------------------------------------
%------------------------------------------------------------------------------------------------------------------------
\subsection{Information theoretic inequalities}
\label{SubSubSec:HellDistFwdL2Dist}
%------------------------------------------------------------------------------------------------------------------------
%------------------------------------------------------------------------------------------------------------------------

	In the following lemma (due to \cite{B2004}) we exploit the boundedness assumption \eqref{Eq:GenUnifBoundG} on $\Gscr$ to show the equivalence between the Hellinger distance appearing in the conclusion of Theorem \ref{Theo:GenHellRates} and the $L^2$-distance on the `regression functions' $\mathscr G(F)$.

\begin{lemma}\label{Lemma:HellAndL2Dist}%----------------------------------------------------------------

	Let the forward map $\Gscr$  satisfy \eqref{Eq:GenUnifBoundG} for some $S>0$. Then, for all $F_1,F_2\in\Fcal$
\begin{equation}
\label{Eq:HellL2}
	\frac{1-e^{-S^2/(2\sigma^2)}}{4S^2}\|\Gscr (F_1)-\Gscr (F_2)\|^2_{L^2(\Dcal)}
	\le
		h^2(p_{F_1},p_{F_2})
	\le
		\frac{1}{4\sigma^2}\|\Gscr (F_1)-\Gscr (F_2)\|_{L^2(\Dcal)}^2.
\end{equation}

\proof%-------------------------------------------------------------------------------------------------------------------

	Note $h^2(p_{F_1},p_{F_2})=2-2\rho(p_{F_1},p_{F_2})$, where
$$
	\rho(p_{F_1},p_{F_2}):=\int_{\R\times\Dcal}
	\sqrt{p_{F_1}(y,x)p_{F_2}(y,x)}dydx
$$
is the Hellinger affinity. Using the expression of the likelihood in \eqref{Eq:SingleLikelihood} (with $\Dcal$ instead of $\Ocal$), the right hand side is seen to be equal to
\begin{align*}
	\int_{\R\times\Dcal}&\frac{1}{\sqrt{2\pi\sigma^2}}
	e^{-\{[y-\Gscr (F_1)(x)]^2-[y-\Gscr (F_2)(x)]^2\}/(4\sigma^2)}dydx\\
%	&=
%		\int_{\R\times\Dcal}\frac{1}{\sqrt{2\pi\sigma^2}}
%		e^{-\{y^2+[\Gscr (F_1)(x)]^2-2y\Gscr (F_1)(x)+y^2
%		+[\Gscr (F_2)(x)]^2-2y\Gscr (F_2)(x)\}/(4\sigma^2)}dydx\\
	&=
		\int_\Dcal e^{-\{[\Gscr (F_1)(x)]^2+[\Gscr (F_2)(x)]^2\}/(4\sigma^2)}
		\Big[\int_\R \frac{e^{-y^2/(2\sigma^2)}}{\sqrt{2\pi\sigma^2}}
		e^{y[\Gscr (F_1)(x)+\Gscr (F_2)(x)]/(2\sigma^2)}dy\Big]dx\\
	&=
		\int_\Dcal e^{-\{[\Gscr (F_1)(x)]^2+[\Gscr (F_2)(x)]^2\}/(4\sigma^2)}
		e^{[\Gscr (F_1)(x)+\Gscr (F_2)(x)]^2/(8\sigma^2)}dx
\end{align*}
having used that the moment generating function of $Z\sim N(0,\sigma^2)$ satisfies $E e^{tZ}=e^{\sigma^2t^2/2}, \ t\in\R$. Thus, the latter integral equals
\begin{align*}
%	\int_\Dcal& e^{-\{2[\Gscr (F_1)(x)]^2+2[\Gscr (F_2)(x)]^2-[\Gscr (F_1)(x)]^2
%	-[\Gscr (F_2)(x)]^2-2\Gscr (F_2)(x)\Gscr (F_2)(x)\}/(8\sigma^2)}dx\\
%	&=
		\int_\Dcal e^{-\{[\Gscr (F_1)(x)]^2+[\Gscr (F_2)(x)]^2
		-2\Gscr (F_2)(x)\Gscr (F_2)(x)\}/(8\sigma^2)}dx
%	&=
%		\int_\Dcal e^{-\{\Gscr (F_1)(x)-\Gscr (F_2)(x)\}^2/(8\sigma^2)}dx\\
	&=
		E^\mu e^{-\{\Gscr (F_1)(X)-\Gscr (F_2)(X)\}^2/(8\sigma^2)}.
\end{align*}

	To derive the second inequality in \eqref{Eq:HellL2}, we use Jensen's inequality to lower bound the expectation in the last line by
\begin{align*}
	e^{-E^\mu\{\Gscr (F_1)(X)-\Gscr (F_2)(X)\}^2/(8\sigma^2)}
	&=
		e^{-\|\Gscr (F_1)-\Gscr (F_2)\|_{L^2(\Dcal)}^2/(8\sigma^2)}.
\end{align*}
Hence
$$
	h^2(p_{F_1},p_{F_2})\le2\Big[1-e^{-\|\Gscr (F_1)-
	\Gscr (F_2)\|_{L^2(\Dcal)}^2/(8\sigma^2)}\Big],
$$
whereby the claim follows using the basic inequality $1-e^{-z/c}\le z/c$, for all $c,z>0$.

	To deduce the first inequality we follow the proof of Proposition 1 in \cite{B2004}: note that for all $0\le z_1< z_2$
% [xxx
% by convexity of $x\in\R\mapsto e^{-x}$
%
$$
	e^{-z_1}
	\le
		\frac{z_1}{z_2}e^{-z_2}+\Big(1-\frac{z_1}{z_2}\Big)
	=
		\frac{e^{-z_2}-1}{z_2}z_1+1.
$$
Then taking $z_1=\{\Gscr (F_1)(X)-\Gscr (F_2)(X)\}^2/(8\sigma^2)$ and $z_2=S^2/(2\sigma^2)$,
\begin{align*}
	E^\mu e^{-\{\Gscr (F_1)(X)-\Gscr (F_2)(X)\}^2/(8\sigma^2)}
	&\le
		\frac{e^{-S^2/(2\sigma^2)}-1}{4S^2}
		\|\Gscr (F_1)-\Gscr (F_2)\|^2_{L^2(\Dcal)}+1
\end{align*}
which in turn yields the result.

\endproof %------------------------------------------------------------------------------------------------------------

\end{lemma}%----------------------------------------------------------------------------------------------------------

	The next lemma bounds the Kullback-Leibler divergences appearing in \eqref{Eq:Bn} in terms of the $L^2$-prediction risk.

\begin{lemma}\label{Lemma:HellingerAndL2SmallBalls}%-------------------------------------------------

	Let the observation $Y_i$ in \eqref{Eq:GenDiscrObs} be generated by some fixed $F_0\in\Fcal$. Then, for each $F\in\Fcal$,
$$
	E^{1}_{F_0}
	\Big[\log\frac{p_{F_0}(Y_1,X_1)}{p_F(Y_1,X_1)}\Big]
	=
		\frac{1}{\sigma^2}\|\Gscr (F_0)-\Gscr (F)\|^2_{L^2(\Dcal)},
$$
and
$$
	E^{1}_{F_0}
	\Big[\log\frac{p_{F_0}(Y_1,X_1)}{p_F(Y_1,X_1)}\Big]^2
	\le
		\frac{2(S^2+\sigma^2)}{\sigma^4}\|\Gscr (F_0)-\Gscr (F)\|^2_{L^2(\Dcal)}.
$$

\proof%------------------------------------------------------------------------------------------------------------------

	If
$
	Y_1=\Gscr (F_0)(X_1)+\sigma W_1,
$
then
\begin{align*}
	\log&\frac{p_{F_0}(Y_1,X_1)}{p_F(Y_1,X_1)}\\
	&=
		-\frac{1}{2\sigma^2}\{[\Gscr (F_0)(X_1)+\sigma W_1-\Gscr (F_0)(X_1)]^2-
		[\Gscr (F_0)(X_1)+\sigma W_1-\Gscr (F)(X_1)]^2\}\\
%	&=
%		-\frac{1}{2\sigma^2}\{(\sigma W_1)^2-
%		[\Gscr (F_0)(X_1)]^2-(\sigma W_1)^2-[\Gscr (F)(X_1)]^2-2\sigma W_1
%		\Gscr (F_0)(X_1)\\
%	&\quad
%		+2 \Gscr (F_0)(X_1)\Gscr (F)(X_1)+ 2\sigma W_1\Gscr (F)(X_1)\}\\
%	&=
%		\frac{1}{2\sigma^2}\{
%		[\Gscr (F_0)(X_1)]^2+[\Gscr (F)(X_1)]^2-2 \Gscr (F_0)(X_1)\Gscr (F)(X_1)\}
%		-\frac{1}{\sigma} W_1\{\Gscr (F_0)(X_1)+\Gscr (F)(X_1)\}
	&=
		\frac{1}{2\sigma^2}\{\Gscr (F_0)(X_1)-\Gscr (F)(X_1)\}^2
		+\frac{1}{\sigma} W_1\{\Gscr (F_0)(X_1)-\Gscr (F)(X_1)\}.
\end{align*}
Hence, since $E W_1=0$ and $X_1\sim\mu$,
\begin{align*}
	E^{1}_{F_0}
	\Big[\log\frac{p_{F_0}(Y_1,X_1)}{p_F(Y_1,X_1)}\Big]
%	&=
%		E^{1}_{F_0}\Big[\frac{1}{2\sigma^2}\{\Gscr (F_0)(X_1)-\Gscr (F)(X_1)\}^2
%		-\frac{1}{\sigma} W_1\{\Gscr (F_0)(X_1)+\Gscr (F)(X_1)\}\Big]\\
	&=
		E^\mu\Big[\frac{1}{2\sigma^2}\{\Gscr (F_0)(X_1)
		-\Gscr (F)(X_1)\}^2\Big]
	=
		\frac{1}{2\sigma^2}\|\Gscr (F_0)-\Gscr (F)\|^2_{L^2(\Dcal)}.
\end{align*}

	On the other hand,
\begin{align*}
	\Big[\log&\frac{p_{F_0}(Y_1,X_1)}{p_F(Y_1,X_1)}\Big]^2\\
	&=
		\Big[\frac{1}{2\sigma^2}\{\Gscr (F_0)(X_1)-\Gscr (F)(X_1)\}^2
		+\frac{1}{\sigma} W_1\{\Gscr (F_0)(X_1)-\Gscr (F)(X_1)\}\Big]^2\\
	&\le
		2\Big[\frac{1}{2\sigma^2}\{\Gscr (F_0)(X_1)-\Gscr (F)(X_1)\}^2\Big]^2+
		2\Big[\frac{1}{\sigma} W_1\{\Gscr (F_0)(X_1)-\Gscr (F)(X_1)\}\Big]^2\\
	% [xxx
	% since $(a+b)^2\le 2a^2+2b^2$ for all $a,b\in\R$.
	% xxx]
	&=
		\frac{2S^2}{\sigma^4}\{\Gscr (F_0)(X_1)-\Gscr (F)(X_1)\}^2+
		\frac{2}{\sigma^2}W_1^2\{\Gscr (F_0)(X_1)-\Gscr (F)(X_1)\}^2,
\end{align*}
whence the second claim follows since $E W_1^2=1$.

\endproof%------------------------------------------------------------------------------------------------------------	
\end{lemma}%----------------------------------------------------------------------------------------------------------

%------------------------------------------------------------------------------------------------------------------------
%------------------------------------------------------------------------------------------------------------------------
\section{Additional background material}
\label{App:BoringThings}
%------------------------------------------------------------------------------------------------------------------------
%------------------------------------------------------------------------------------------------------------------------

In this final appendix we collect some standard materials used in the proofs for convenience of the reader.

\begin{example}\label{Ex:LinkFunction}\normalfont%-------------------------------------------------------

Take
$$
	\phi:\R\to(0,\infty), \quad \phi(t)=\frac{1}{1-t}1_{\{t<0\}}+(1+t)1_{\{t\ge0\}}, 
$$
and let $\psi:\R\to[0,\infty)$ be a smooth compactly supported function such that $			\int_\R\psi(t)dt=1$. Define for any $K_{min}\in(0,1)$
\begin{equation}
\label{Eq:Phi}
	\Phi(t)=K_{min}+\frac{1-K_{min}}{\psi * \phi(0)}\psi * \phi (t), \quad t\in\R.
\end{equation}
Then it is elementary to check that $\Phi$ is a regular link function that satisfies Condition \ref{Cond:LinkFunction2} (with $a=2$). 
\end{example}%-------------------------------------------------------------------------------------------------------

\begin{example}\label{Ex:CuttedMaternProcess} \normalfont%-----------------------------------------

	For any real $\alpha>d/2$, the Whittle-Matérn process with index set $\Ocal$ and regularity $\alpha-d/2>0$ (cf. Example 11.8 in \cite{GvdV17})  is the stationary centred Gaussian process $M=\{M(x), \ x\in\Ocal\}$ with covariance kernel 
$$
	K(x,y)=\int_{\R^d}e^{-i\langle x-y,\xi\rangle_{\R^d}}\mu(d\xi), 
	\quad \mu(d\xi)=(1+\|\xi\|^2_{\R^d})^{-\alpha}d\xi,~~x,y \in \Ocal.
$$
From the results in Chapter 11 in \cite{GvdV17} we see that the RKHS of $(M(x): x \in \Ocal)$ equals the set of restrictions to $\Ocal$ of elements in the Sobolev space $H^\alpha(\mathbb R^d)$, which equals, with equivalent norms, the space $H^\alpha(\Ocal)$ (since $\Ocal$ has a smooth boundary). Moreover, Lemma I.4 in \cite{GvdV17} shows that $M$ has a version with paths belonging almost surely to $C^{\beta'}$ for all $\beta'<\alpha-d/2$. Let now $K\subset\Ocal$ be a nonempty compact set, and let $M$ be a $C^{\beta'}$-smooth version of a Whittle-Matérn process on $\Ocal$ with RKHS $H^\alpha(\Ocal)$. Taking $F'=\chi M$ implies (cf. Exercise 2.6.5 in \cite{GN16}) that $\Pi'=\Lcal(F')$ defines a centred Gaussian probability measure supported on $C^{\beta'}$, whose RKHS is given by
$$\Hcal=\{\chi F,\ F\in H^\alpha(\Ocal)\},$$
and the RKHS norm satisfies that for all $F \in H^\alpha(\Ocal)$ there exists $F^* \in H^\alpha(\Ocal)$ such that $\chi F = \chi F^*$ and
$$\|\chi F\|_\Hcal = \|F^*\|_{H^\alpha(\Ocal)}.$$
Thus if $F' = \chi F$ is an arbitrary element of $\Hcal$, then 
$$\|F'\|_{H^\alpha}=\| \chi F^*\|_{H^\alpha}\lesssim \|F^*\|_{H^\alpha}=\|F'\|_{\Hcal},$$
which shows that $\Hcal$ is continuously embedded into $H^\alpha_c(\Ocal)$. 
\end{example}%-------------------------------------------------------------------------------------------------------

\begin{remark}\label{Rem:WaveletSeries}\normalfont%---------------------------------------------------

	Let $\{\Psi_{\ell r}, \ \ell\ge-1, r\in\Z^d\}$ be an orthonormal basis of $L^2(\R^d)$ composed of $S$-regular and compactly supported Daubechies wavelets (see  Chapter 4 in \cite{GN16} for construction and properties). For each $0\le \alpha\le S$, we have
$$
	H^\alpha(\R^d)
	=
		\Big\{F\in L^2(\R^d):
		\sum_{\ell, r}2^{2\ell \alpha}\langle F,\Psi_{\ell r}\rangle^2_{L^2(\R^d)}<\infty\Big\},
$$
and the square root of the latter series defines an equivalent norm to $\|\cdot\|_{H^\alpha(\R^d)}$. Note that $S>0$ can be taken arbitrarily large.

For any $\alpha\ge0$ the Gaussian random series
$$
	\bar F_j
	=
		\sum_{\ell \le j, r\in\Rcal_\ell} F_{\ell r}2^{-\ell \alpha}\Psi_{\ell r},
	 \quad
	  F_{\ell r}\iid N(0,1)
$$
defines a centred Gaussian probability measure supported on the finite-dimensional space $\bar\Hcal_j$  spanned by the $\{\Psi_{\ell r}, \ell \le j, r \in \Rcal_\ell\}$, and its RKHS equals $\bar\Hcal_j$ endowed with norm
$$
	\| \bar H_j \|_{\bar\Hcal_j}^2
	=
		\sum_{\ell \le j, r\in\Rcal_\ell} 2^{2\ell \alpha}H^2_{\ell r}
	=
		\| \bar H_j\|^2_{H^\alpha(\R^d)}
	\quad
	\forall \bar H_j\in\bar\Hcal_j
$$ 
(cf. Example 2.6.15 in \cite{GN16}). Basic wavelet theory implies $\textnormal{dim}(\bar\Hcal_j)\lesssim 2^{jd}$.

	If we now fix compact $K'\subset \Ocal$ such that $K\subsetneq K'$, and consider a cut-off function $\chi\in C^\infty_c(\Ocal)$ such that $\chi=1$ on $K'$, then multiplication by $\chi$ is a bounded linear operator
$
	\chi : H^s(\R^d)\to H^s_c(\Ocal).
$
It follows that the random function
$$
	F_j
	=
		\chi (\bar F_j)
	=
		\sum_{\ell \le j, r\in\Rcal_\ell} F_{\ell r}2^{-\ell \alpha}\chi\Psi_{\ell r},
	 \quad
	  F_{\ell r}\iid N(0,1)
$$
defines, according to Exercise 2.6.5 in \cite{GN16}, a centred Gaussian probability measure $\Pi_j=\Lcal(F_j)$ supported on the finite dimensional subspace $\Hcal_j$ from (\ref{Eq:SubspaceHj}), with RKHS norm satisfying
\begin{align}
\label{Eq:FiniteRKHSNorm}
	\Big\|\chi\Big(\sum_{\ell \le j, r\in\Rcal_\ell} H_{\ell r}\Psi_{\ell r}\Big)\Big\|_{\Hcal_j}
	&\le
		\Big\|\sum_{\ell \le j, r\in\Rcal_\ell} H_{\ell r}\Psi_{\ell r}\|_{\bar\Hcal_j}
%	&=
%		\Big\|\sum_{\ell \le j, r\in\Rcal_\ell} H_{\ell r}\Psi_{\ell r}\Big\|_{H^\alpha(\R^d)}
	=
		\sqrt{\sum_{\ell \le j, r\in\Rcal_\ell} 2^{2\ell\alpha}H^2_{\ell r}}.
\end{align}
Arguing as in the previous remark one shows further that for some constant $c>0$,
\begin{equation}
\label{Eq:Norms}
	\|H_j\|_{H^\alpha(\Ocal)}\le c\|H_j\|_{\Hcal_j} \quad \forall H_j\in\Hcal_j.
\end{equation}
%Indeed, for $H_j\in \Hcal_j$, there exists $\bar H_j\in \bar\Hcal_j$ such that $H_j=\chi \bar H_j$ and $\|H_j\|_{\Hcal_j}=\|\bar \Hcal_j\|_{\bar \Hcal_j}=\|\bar \Hcal_j\|_{H^\alpha(\R^d)}$. Then, in view of the continuity of $\chi$
%$$
%	\|H_j\|_{H^\alpha(\Ocal)}
%	=
%		\|\chi(\bar H_j)\|_{H^\alpha(\Ocal)}
%	\lesssim
%		\|\bar H_j\|_{H^\alpha(\R^d)}
%	=
%		\|H_j\|_{\Hcal_j}.
%$$

\end{remark}%---------------------------------------------------------------------------------------------------------

\begin{remark}\label{Rem:ApproxProp}\normalfont%-------------------------------------------------------

	Using the notation of the previous remark, for fixed $F_0\in H^\alpha_K(\Ocal)$, consider the finite-dimensional approximations 
\begin{equation}
\label{Eq:Projections}
	P_{\Hcal_j}(F_0)
	=
		\sum_{\ell\le j, r\in\Rcal_\ell}\langle F_0,\Psi_{\ell r}\rangle_{L^2}
		\chi\Psi_{\ell r}
	\in \Hcal_j,
	\quad
	j\in\N.
\end{equation}
Then in view of \eqref{Eq:FiniteRKHSNorm}, we readily check that for each $j\ge1$
\begin{equation}
\label{Eq:AnotherBound}
	\|P_{\Hcal_j}(F_0)\|_{\Hcal_j}
	\le
		\sqrt{\sum_{\ell\le j, r\in\Rcal_\ell} 2^{2\ell\alpha}\langle F_0,\Psi_{\ell r}\rangle_{L^2}^2}
	\le
		\|F_0\|_{H^\alpha(\Ocal)}
	<\infty.
\end{equation}
Also, for each $\kappa\ge0$, and any $G\in H^\kappa(\Ocal)$, we see that (implicitly extending to $0$ on $\R^d\backslash\Ocal$ functions that are compactly supported inside $\Ocal$)
$$
	\langle F_0-P_{\Hcal_j}(F_0),G\rangle_{L^2(\Ocal)}
	=
		\langle F_0-P_{\Hcal_j}(F_0),\chi'G\rangle_{L^2(\R^d)}
$$
where $\chi'\in C^\infty_c(\Ocal)$, with $\chi'=1$ on $\textnormal{supp}(\chi)$. 
We also note that, in view of the localisation properties of Daubechies wavelets, for some $J_{min}\in\N$ large enough, if $\ell\ge J_{min}$ and the support of $\Psi_{\ell r}$ intersects $K$, then necessarily $\textnormal{supp}(\Psi_{\ell r})\subseteq K'$, so that
$$
	\chi\Psi_{\ell r}=\Psi_{\ell r} \quad \forall \ell \ge J_{min}, \ r\in\Rcal_\ell.
$$
Therefore, for $j\ge J_{min}$, by Parseval's identity and the Cauchy-Schwarz inequality
\begin{align*}
	\langle F_0&-P_{\Hcal_j}(F_0),\chi'G\rangle_{L^2(\R^d)}\\
	&=
		\sum_{\ell'>j, r'\in\Rcal_\ell}2^{\ell \alpha}\langle F_0,\Psi_{\ell' r'}\rangle_{L^2(\R^d)}
		2^{\ell'\kappa}\langle \chi' G,\Psi_{\ell' r'}\rangle_{L^2(\R^d)}2^{-\ell'(\alpha+\kappa)}
		\\
	&\le
		2^{-j(\alpha+\kappa)}\sqrt{\sum_{\ell'>j, r'\in\Rcal_\ell}2^{2\ell \alpha}
		\langle F_0,\Psi_{\ell' r'}\rangle_{L^2(\R^d)}^2}\sqrt{\sum_{\ell'>j, r'\in\Rcal_\ell}
		2^{2\ell \kappa}\langle \chi' G,\Psi_{\ell' r'}\rangle_{L^2(\R^d)}^2}\\
	&\le
		2^{-j(\alpha+\kappa)}\|F_0\|_{H^\alpha(\Ocal)}\|\chi' G\|_{H^\kappa(\R^d)}.
 \end{align*}
 It follows by duality that for all $j$ large enough
 \begin{equation}
 \label{Eq:DualNormApprox}
 	\|F_0-P_{\Hcal_j}( F_0)\|_{(H^\kappa(\Ocal))^*}
%	=
%		\sup_{G\in H^\kappa(\Ocal), \ \|G\|_{H^\kappa(\Ocal)}\le1}|\langle\bar F-P_J\bar 		F,G\rangle_{L^2(\Ocal)}|
	\lesssim
		2^{-j(\alpha+\kappa)}\|F_0\|_{H^\alpha(\Ocal)}.
\end{equation}
We conclude remarking that
\begin{equation}
\label{Eq:EquivOfNorms}
	\|F\|_{H^\alpha(\Ocal)}\lesssim 2^{j\alpha}\|F\|_{L^2(\Ocal)}, \quad
	\forall F\in \Hcal_j,\ j\ge J_{min}.
\end{equation}
Indeed, let $j\ge J_{min}$, and fix $F\in\Hcal_j$; then
\begin{align*}
	F
	&=
		P_{\Hcal_{J_{min}}} (F)+(F-P_{\Hcal_{J_{min}}}(F))
	=
		\sum_{\ell\le J_{min},r\in\Rcal_\ell}F_{\ell r}\chi\Psi_{\ell r} 
		+ \sum_{J_{min}<\ell\le j,r\in\Rcal_\ell}F_{\ell r}\Psi_{\ell r}.
\end{align*}
But as $\Hcal_{J_{min}}$ is a fixed finite dimensional subspace,  then we have $\|P_{\Hcal_{J_{min}}}(F)\|_{H^s(\Ocal)}\lesssim\|P_{\Hcal_{J_{min}}} (F)\|_{L^2(\Ocal)}\le \|F\|_{L^2(\Ocal)}$ for some fixed multiplicative constant only depending on $J_{min}$. On the other hand, we also have
\begin{align*}
	\|F-P_{\Hcal_{J_{min}}}( F)\|^2_{H^\alpha(\Ocal)}
	&=
		\sum_{J_{min}<\ell\le j,r\in\Rcal_\ell} 2^{2\ell \alpha}F_{\ell r}^2\\
%	&\le
%		2^{2J s} \sum_{\ell',r'}\Big\langle\sum_{J_{min}<\ell\le j,r\in\Rcal_\ell} F_{\ell r}^2 
%		\Psi_{\ell r},\Psi_{\ell' r'}\Big\rangle_{L^2(\R^d)}^2\\
	&\le
		2^{2j \alpha}\|F-P_{\Hcal_{J_{min}}}( F)\|^2_{L^2(\Ocal)}\le2^{2j \alpha}\|F\|^2_{L^2(\Ocal)},
\end{align*}
yielding \eqref{Eq:EquivOfNorms}.

\end{remark}%---------------------------------------------------------------------------------------------------------

\begin{example}\label{poissond}\normalfont%----------------------------------------------------------------
	Consider the integer-valued random variable
$$
	J=\lfloor \log_2(\phi^{-1}(T)^{1/d})\rfloor+1, \quad T\sim\textnormal{Exp}(1),
$$
where $\phi(x)=x\log x, \ x\ge1$. Then for any $j\ge1$
 \begin{align*}
	\Pr(J> j)
%	&=
%		\Pr(\lfloor \log_2(\phi^{-1}(T)^{1/d}) \rfloor> j-1)\\
%	&=
%		\Pr( \log_2(\phi^{-1}(T)^{1/d})\ge j)\\
%	&=
%		\Pr( \phi^{-1}(T)^{1/d}\ge 2^j)\\
	&=
		\Pr( \phi^{-1}(T)\ge 2^{jd}) =
		\Pr( T\ge 2^{jd}\log 2^{jd}) =
		e^{-2^{jd}\log 2^{jd}}.
\end{align*}
On the other hand, since $e^{-2^{jd}(1-2^{-d})\log2^{(j-1)d}}\to0$ as $j\to\infty$,
\begin{align*}
	\Pr(J= j) 	&=
		\Pr( 2^{(j-1)d}\le \phi^{-1}(T)< 2^{jd})=
		e^{-2^{(j-1)d}\log2^{(j-1)d}}+1-e^{-2^{jd}\log 2^{jd}}-1\\
%	&=
%		e^{-2^{(j-1)d}\log 2^{(j-1)d}}
%		\Bigg(1-\frac{e^{2^{(j-1)d}\log2^{(j-1)d}}}{e^{2^{jd}\log 2^{jd}}}\Bigg)\\
	&\ge
%		e^{-2^{(j-1)d}\log 2^{(j-1)d}}
%		\Bigg(1-\frac{e^{2^{(j-1)d}\log2^{(j-1)d}}}{e^{2^{jd}\log 2^{(j-1)d}}}\Bigg)\\
%%	&=
%%		e^{-2^{(j-1)d}\log 2^{(j-1)d}}
%%		(1-e^{[2^{(j-1)d}-2^{jd}]\log2^{(j-1)d}})\\
%	&=
		e^{-2^{(j-1)d}\log 2^{(j-1)d}}
		(1-e^{-2^{jd}(1-2^{-d})\log2^{(j-1)d}}) \gtrsim e^{-2^{jd}\log 2^{jd}}.
\end{align*}
\end{example}%-------------------------------------------------------------------------------------------------------

\smallskip

\textbf{Acknowledgement.} The authors are grateful to three anonymous referees for critical remarks and suggestions, as well as to Sven Wang for helpful discussions. We  further acknowledge support by the European Research Council under ERC grant agreement No.647812 (UQMSI).

\bibliographystyle{acm}

\bibliography{GNWreferences}

\begin{thebibliography}{10}

\bibitem{AN19}
{\sc Abraham, K., and Nickl, R.}
\newblock On statistical {C}ald\'eron problems.
\newblock {\em arXiv:1906.03486\/} (2019).

\bibitem{S13}
{\sc Agapiou, S., Larsson, S., and Stuart, A.~M.}
\newblock Posterior contraction rates for the {B}ayesian approach to linear
  ill-posed inverse problems.
\newblock {\em Stochastic Process. Appl. 123}, 10 (2013), 3828--3860.

\bibitem{A86}
{\sc Alessandrini, G.}
\newblock An identification problem for an elliptic equation in two variables.
\newblock {\em Annali di Matematica Pura ed Applicata 145}, 1 (1986), 265--295.

\bibitem{AMOS19}
{\sc Arridge, S., Maass, P., \"{O}ktem, O., and Sch\"{o}nlieb, C.-B.}
\newblock Solving inverse problems using data-driven models.
\newblock {\em Acta Numer. 28\/} (2019), 1--174.

\bibitem{BB18}
{\sc Benning, M., and Burger, M.}
\newblock Modern regularization methods for inverse problems.
\newblock {\em Acta Numerica 27\/} (2018), 1--111.

\bibitem{BGLFS17}
{\sc Beskos, A., Girolami, M., Lan, S., Farrell, P., and Stuart, A.}
\newblock Geometric {MCMC} for infinite-dimensional inverse problems.
\newblock {\em Journal of Computational Physics 335\/} (2017), 327--351.

\bibitem{B2004}
{\sc Birg{\'e}, L.}
\newblock Model selection for gaussian regression with random design.
\newblock {\em Bernoulli 10\/} (2004), 1039--1051.

\bibitem{BHM04}
{\sc Bissantz, N., Hohage, T., and Munk, A.}
\newblock Consistency and rates of convergence of nonlinear {T}ikhonov
  regularization with random noise.
\newblock {\em Inverse Problems 20}, 6 (2004), 1773--1789.

\bibitem{BHMR07}
{\sc Bissantz, N., Hohage, T., Munk, A., and Ruymgaart, F.}
\newblock Convergence rates of general regularization methods for statistical
  inverse problems and applications.
\newblock {\em SIAM Journal on Numerical Analysis 45}, 6 (2007), 2610--2636.

\bibitem{BCDPW17}
{\sc Bonito, A., Cohen, A., DeVore, R., Petrova, G., and Welper, G.}
\newblock Diffusion coefficients estimation for elliptic partial differential
  equations.
\newblock {\em SIAM J. Math. Anal. 49}, 2 (2017), 1570--1592.

\bibitem{G19}
{\sc Briol, F.-X., Oates, C.~J., Girolami, M., Osborne, M.~A., and Sejdinovic,
  D.}
\newblock Probabilistic integration: a role in statistical computation?
\newblock {\em Statist. Sci. 34}, 1 (2019), 1--22.

\bibitem{CMPS16}
{\sc Conrad, P.~R., Marzouk, Y.~M., Pillai, N.~S., and Smith, A.}
\newblock Accelerating asymptotically exact {MCMC} for computationally
  intensive models via local approximations.
\newblock {\em Journal of the American Statistical Association 111}, 516
  (2016), 1591--1607.

\bibitem{CRSW13}
{\sc Cotter, S., Roberts, G., Stuart, A., and White, D.}
\newblock {MCMC} methods for functions: Modifying old algorithms to make them
  faster.
\newblock {\em Statistical Science 28}, 3 (2013), 424--446.

\bibitem{DS11}
{\sc Dashti, M., and Stuart, A.~M.}
\newblock Uncertainty quantification and weak approximation of an elliptic
  inverse problem.
\newblock {\em SIAM J. Numer. Anal. 49}, 6 (2011), 2524--2542.

\bibitem{DS16}
{\sc Dashti, M., and Stuart, A.~M.}
\newblock The {B}ayesian approach to inverse problems.
\newblock {\em In: Handbook of Uncertainty Quantification, Editors R. Ghanem,
  D. Higdon and H. Owhadi, Springer\/} (2016).

\bibitem{D88}
{\sc Diaconis, P.}
\newblock {Bayesian Numerical Analysis}.
\newblock In {\em Statistical Decision Theory and Related Topics IV\/} (1988),
  J.~Berger and S.~Gupta, Eds., Springer, New York, pp.~163--175.

\bibitem{EHN96}
{\sc Engl, H.~W., Hanke, M., and Neubauer, A.}
\newblock {\em Regularization of Inverse Problems}.
\newblock Kluwer Academic Publishers Group, 1996.

\bibitem{F83}
{\sc Falk, R.~S.}
\newblock Error estimates for the numerical identification of a variable
  coefficient.
\newblock {\em Mathematics of Computation 40}, 162 (1983), 537--546.

\bibitem{GvdV17}
{\sc Ghosal, S., and van~der Vaart, A.~W.}
\newblock {\em Fundamentals of Nonparametric Bayesian Inference}.
\newblock Cambridge University Press, New York, 2017.

\bibitem{GT98}
{\sc Gilbarg, D., and Trudinger, N.~S.}
\newblock {\em Elliptic partial differential equations of second order}.
\newblock Springer-Verlag, Berlin-New York, 1998.

\bibitem{GN16}
{\sc Gin\'e, E., and Nickl, R.}
\newblock {\em Mathematical foundations of infinite-dimensional statistical
  models}.
\newblock Cambridge University Press, New York, 2016.

\bibitem{HSV14}
{\sc Hairer, M., Stuart, A., and Vollmer, S.}
\newblock Spectral gaps for a {M}etropolis-{H}astings algorithm in infinite
  dimensions.
\newblock {\em The Annals of Applied Probability 24}, 6 (2014), 2455--2490.

\bibitem{HNS95}
{\sc Hanke, M., Neubauer, A., and Scherzer, O.}
\newblock A convergence analysis of the {L}andweber iteration for nonlinear
  ill-posed problems.
\newblock {\em Numer. Math. 72}, 1 (1995), 21--37.

\bibitem{HS85}
{\sc Hoffmann, K.~H., and Sprekels, J.}
\newblock On the identification of coefficients of elliptic problems by
  asymptotic regularization.
\newblock {\em Numerical Functional Analysis and Optimization 7}, 2-3 (1985),
  157--177.

\bibitem{HP08}
{\sc Hohage, T., and Pricop, M.}
\newblock Nonlinear {T}ikhonov regularization in {H}ilbert scales for inverse
  boundary value problems with random noise.
\newblock {\em Inverse Probl. Imaging 2}, 2 (2008), 271--290.

\bibitem{IK94}
{\sc Ito, K., and Kunisch, K.}
\newblock On the injectivity and linearization of the coefficient-to-solution
  mapping for elliptic boundary value problems.
\newblock {\em Journal of Mathematical Analysis and Applications 188}, 3 (dec
  1994), 1040--1066.

\bibitem{KS04}
{\sc Kaipio, J., and Somersalo, E.}
\newblock {\em Statistical and Computational Inverse Problems}.
\newblock No.~160 in Applied Mathematical Sciences. Springer-Verlag New York,
  2004.

\bibitem{KNS08}
{\sc Kaltenbacher, B., Neubauer, A., and Scherzer, O.}
\newblock In {\em Iterative Regularization Methods for Nonlinear Ill-Posed
  Problems\/} (2008), Radon Series on Computational and Applied Mathematics.

\bibitem{KSS09}
{\sc Kaltenbacher, B., Sch\"{o}pfer, F., and Schuster, T.}
\newblock Iterative methods for nonlinear ill-posed problems in {B}anach
  spaces: convergence and applications to parameter identification problems.
\newblock {\em Inverse Problems 25\/} (2009).

\bibitem{KLS16}
{\sc Kekkonen, H., Lassas, M., and Siltanen, S.}
\newblock Posterior consistency and convergence rates for {B}ayesian inversion
  with hypoelliptic operators.
\newblock {\em Inverse Problems 32}, 8 (2016), 085005, 31.

\bibitem{K2015}
{\sc Knapik, B., Szab{\`o}, B., van~der Vaart, A.~W., and van Zanten, H.}
\newblock Bayes procedures for adaptive inference in inverse problems for the
  white noise model.
\newblock {\em Probab. Theory Relat. Fields}, 164 (2015), 771--813.

\bibitem{K11}
{\sc Knapik, B., van~der Vaart, A.~W., and van Zanten, J.~H.}
\newblock Bayesian inverse problems with {G}aussian priors.
\newblock {\em Ann. Statist. 39}, 5 (2011), 2626--2657.

\bibitem{K01}
{\sc Knowles, I.}
\newblock Parameter identification for elliptic problems.
\newblock {\em Journal of Computational and Applied Mathematics 131}, 1 (2001),
  175 -- 194.

\bibitem{KL88}
{\sc Kohn, R.~V., and Lowe, B.~D.}
\newblock A variational method for parameter identification.
\newblock {\em ESAIM: Mathematical Modelling and Numerical Analysis -
  Mod\'elisation Math\'ematique et Analyse Num\'erique 22}, 1 (1988), 119--158.

\bibitem{KS85}
{\sc Kravaris, C., and Seinfeld, J.}
\newblock Identification of parameters in distributed parameter systems by
  regularization.
\newblock {\em SIAM Journal on Control and Optimization 23\/} (03 1985).

\bibitem{LL99}
{\sc Li, W.~V., and Linde, W.}
\newblock Approximation, metric entropy and small ball estimates for {G}aussian
  measures.
\newblock {\em Ann. Probab. 27}, 3 (1999), 1556--1578.

\bibitem{Lions1972}
{\sc Lions, J.~L., and Magenes, E.}
\newblock {\em Non-Homogeneous Boundary Value Problems and Applications}, 1~ed.
\newblock Grundlehren der mathematischen Wissenschaften. Springer-Verlag Berlin
  Heidelberg, 1972.

\bibitem{MNP19b}
{\sc Monard, F., Nickl, R., and Paternain, G.~P.}
\newblock Consistent inversion of noisy non-abelian {X}-ray transforms.
\newblock {\em arXiv:1905.00860\/} (2019).

\bibitem{MNP17}
{\sc Monard, F., Nickl, R., and Paternain, G.~P.}
\newblock Efficient nonparametric {B}ayesian inference for {$X$}-ray
  transforms.
\newblock {\em Ann. Statist. 47}, 2 (2019), 1113--1147.

\bibitem{N17}
{\sc Nickl, R.}
\newblock {Bernstein}-von {Mises} theorems for statistical inverse problems
  {I}: {Schr{\"o}dinger} equation.
\newblock {\em Journal of the European Mathematical Society (JEMS), to
  appear\/} (2018).

\bibitem{NS17}
{\sc Nickl, R., and S\"ohl, J.}
\newblock Nonparametric {B}ayesian posterior contraction rates for discretely
  observed scalar diffusions.
\newblock {\em Ann. Statist. 45}, 4 (2017), 1664--1693.

\bibitem{NS19}
{\sc Nickl, R., and S\"{o}hl, J.}
\newblock Bernstein--von {M}ises theorems for statistical inverse problems
  {II}: compound {P}oisson processes.
\newblock {\em Electron. J. Stat. 13}, 2 (2019), 3513--3571.

\bibitem{NVW18}
{\sc Nickl, R., van~de Geer, S., and Wang, S.}
\newblock Convergence rates for penalised least squares estimators in
  {PDE}-constrained regression problems.
\newblock {\em SIAM/ASA Journal of Uncertainty Quantification\/} (to appear).

\bibitem{Q00}
{\sc Qi-nian, J.}
\newblock Error estimates of some newton-type methods for solving nonlinear
  inverse problems in {H}ilbert scales.
\newblock {\em Inverse Problems 16}, 1 (2000), 187--197.

\bibitem{R13}
{\sc Ray, K.}
\newblock Bayesian inverse problems with non-conjugate priors.
\newblock {\em Electron. J. Stat. 7\/} (2013), 2516--2549.

\bibitem{R81}
{\sc Richter, G.~R.}
\newblock An inverse problem for the steady state diffusion equation.
\newblock {\em SIAM J. Appl. Math. 41}, 2 (1981), 210--221.

\bibitem{SS12}
{\sc Schwab, C., and Stuart, A.~M.}
\newblock Sparse deterministic approximation of {B}ayesian inverse problems.
\newblock {\em Inverse Problems 28}, 4 (2012), 045003, 32.

\bibitem{S10}
{\sc Stuart, A.~M.}
\newblock Inverse problems: a {B}ayesian perspective.
\newblock {\em Acta Numer. 19\/} (2010), 451--559.

\bibitem{Triebel1978}
{\sc Triebel, H.}
\newblock {\em Interpolation Theory, Function Spaces, Differential Operators}.
\newblock North-Holland Publishing Company, 1978.

\bibitem{vdVvZ08}
{\sc van~der Vaart, A.~W., and van Zanten, J.~H.}
\newblock Rates of contraction of posterior distributions based on {G}aussian
  process priors.
\newblock {\em Ann. Statist. 36}, 3 (2008), 1435--1463.

\bibitem{vdVvZ09}
{\sc van~der Vaart, A.~W., and van Zanten, J.~H.}
\newblock Adaptive {B}ayesian estimation using a {G}aussian random field with
  inverse gamma bandwidth.
\newblock {\em Ann. Statist. 37}, 5B (2009), 2655--2675.

\bibitem{V13}
{\sc Vollmer, S.~J.}
\newblock Posterior consistency for {B}ayesian inverse problems through
  stability and regression results.
\newblock {\em Inverse Problems 29}, 12 (2013), 125011, 32.

\end{thebibliography}

\end{document}